\documentclass[11pt]{article}
\usepackage{lscape}
 \usepackage[ansinew]{inputenc}
 \usepackage[dvips]{graphicx}
\usepackage[psamsfonts]{amssymb} \usepackage[all]{xy}  \usepackage[all]{xy}
 \usepackage[psamsfonts]{eucal}
 \usepackage[psamsfonts]{eucal}
\usepackage{amssymb, amsmath}
\usepackage{amsthm}
\usepackage{geometry}
\usepackage{enumerate}
\usepackage{float}
\newcommand{\be}{\begin{equation}}
\newcommand{\ee}{\end{equation}}

\newcommand{\bd}{\begin{displaymath}}
\newcommand{\ed}{\end{displaymath}}
\newcommand{\ben}{\begin{enumerate}}
\newcommand{\een}{\end{enumerate}}

\newtheorem{Prop}{Proposition}

\title{Sullivan completions}
\author{Y. F\'elix and S. Halperin}
\begin{document}
\maketitle
 
 \begin{abstract} The Sullivan construction   associates to each path connected space or connected simplicial set, $X$, a special cdga, its minimal model $(\land V,d)$, and to each such cdga $\land W$ its geometric realisation $\langle \land W\rangle$. The composite of these constructions is  the Sullivan completion, $X_{\mathbb Q}$, of $X$. In this paper we give a survey of  the main   properties   of  Sullivan completions, and include   explicit examples.    \end{abstract}
 
 \section{Introduction} 
 The fundamental contributions of Quillen \cite{Q} and Sullivan \cite{nSu} in 1970 were initial landmarks in formalizing the rationalization of a path connected space. In fact, Quillen completely solved the problem for a simply connected space, identifying the rationalization as the realization of a differential graded Lie algebra. Sullivan then formalized his approach via his minimal models, \cite{Su}.

 Intuitively, the homology and homotopy groups of the rationalization of a space should be the tensor product of the original groups with $\mathbb Q$, and this is the case with Quillen's construction. Since this is not possible for non-simply connected spaces, $X$, one looks instead for useful "completions", inspired in particular by the classical completion $B\to \widehat{B} $ of an augmented algebra $B\stackrel{\varepsilon}{\longrightarrow} \mathbb Q$:
 $$\widehat{B}:= \varprojlim_k B/I^k,$$
 $I^k$ denoting the $k^{th}$ power of the augmentation ideal $I:=\mbox{ker}\,\varepsilon$.
 
For path  connected spaces, or connected simplicial sets,  
  the $R$-completions $R_\infty(X)$ of Bousfield-Kan \cite{BK}, defined for any commutative ring $R$, are prime examples, and it follows from (\cite[(4.3), p.137]{BK}) that up to homotopy these coincide with Quillen's construction for simply connected spaces when $R= \mathbb Q$. On the other hand, there is a similarity to the constructions of $\mathbb Q_\infty (X)$ and  of Sullivan completions $X_{\mathbb Q}$. For instance, Sullivan completions can be described as a functor $\ell_{\mathbb Q} : X\to X_{\mathbb Q}$, characterized by a universal property analogous to that of $\mathbb Q_\infty(X)$ established in (\cite[III (6.2)]{BK}). 
  The Bousfield-Kan completion of $X$ is the inverse limit
 $$X\to R_\infty(X) = \varprojlim_kR_k(X)$$
 of a tower of fibrations, and inducing an isomorphism $\varinjlim_k H^*(R_k(X);R) \stackrel{\cong}{\to} H^*(X;R)$. Analogously, the Sullivan completion of $X$ is the inverse limit
 $$\ell: X\to X_{\mathbb Q}= \varprojlim_\alpha X_\alpha$$
 of a specific inverse system inducing an isomorphism 
 $\varinjlim_\alpha H^*(X_\alpha;\mathbb Q)\stackrel{\cong}{\to} H^*(X;\mathbb Q)$. Here the $X_\alpha$ are nilpotent spaces defined via the Sullivan model of $X$, and for which $\oplus_{i\geq 2}\pi_i(X)$ and each $\pi_1(X)^{n}/ \pi_1(X)^{n+1}$  is a  finite dimensional   rational vector space. Here $\pi_1(X)^{n}$ denotes the $n^{th}$ term of the lower central series of $\pi_1(X)$.

Finally, as Bousfield has shown (\cite{2a}), the completions $\mathbb Q_\infty(X)$ and $X_{\mathbb Q}$ are directly related: the map $X\to X_{\mathbb Q}$ factors up to homotopy 	as $X\to {\mathbb Q}_\infty (X) \to X_{\mathbb Q}$. Moreover, if in the terminology of \cite{BK} $X$ is $\mathbb Q$-good, then $\mathbb Q_{\infty}(X)_{\mathbb Q}\to X_{\mathbb Q}$ is a homotopy equivalence. Furthermore, (\cite[Theorem 12.2]{BG})    if $H_*(X;\mathbb Q)$ is a graded vector space of finite type then $\mathbb Q_\infty(X)\to X_{\mathbb Q}$ is a homotopy equivalence. 
   Finally we note that for a map $f : X\to Y$ between path connected spaces the following conditions are equivalent: (i) $H^*(f;\mathbb Q)$ is an isomorphism, (ii) $\mathbb Q_\infty (f) : \mathbb Q_\infty (X)\to \mathbb Q_\infty (Y)$ is a homotopy equivalence, and (iii) $f_{\mathbb Q}: X_{\mathbb Q}\to Y_{\mathbb Q}$ is a homotopy equivalence.

 Here we shall focus on the properties of Sullivan completions, $X_{\mathbb Q}$. These are constructed as the simplicial realization of a minimal Sullivan model $(\land V,d)$ for $X$. Because $H(\land V,d)\cong H^*(X;\mathbb Q)$ this model provides an often computable approach to the cohomology of $X$.
 On the other hand, 
 there is a  natural bijection
 $$\pi_*(X_{\mathbb Q})\cong \mbox{Hom}(V,\mathbb Q)$$
which provides access to the structure and properties of $\pi_*(X_{\mathbb Q})$.

 Some Sullivan completions are very simple. For instance $S^1_{\mathbb Q}$ is the rationalization of $S^1$,  constructed by a telescope process. But in general, Sullivan completions are more complicated and with quite mysterious cohomology. For example, Ivanov and Mikhailov have recently proved \cite{IM} that $H_2((S^1\vee S^1)_{\mathbb Q};\mathbb Q)$ is uncountable. Moreover it is unclear whether $H_{\geq 1}(X_{\mathbb Q};\mathbb Z)$ is a rational vector space, though this is true for the $X_\alpha$.

For simplicity we adopt the following

\vspace{3mm}\noindent {\bf Convention} \begin{enumerate}
\item[(i)] By "space" we mean either a CW complex or a simplicial set.
\item[(ii)] For a space   $X$ we write
$$H(X):= H^*(X;\mathbb Q).$$
\item[(iii)] We write
$$(-)^\vee := \mbox{Hom}(-, \mathbb Q).$$
\item[(iv)] Where there is no ambiguity we suppress the differentials from the notation for a complex and write $A$ instead of $(A,d)$.

\end{enumerate}
 
 \vspace{2mm}Finally, our thanks to Pete Bousfield, \cite{2a}, for a number of helpful suggestions, including the observation that $X\to X_{\mathbb Q}$ factors through ${\mathbb Q}_\infty (X)$.

 \section{Basic constructions} 
 
 We   briefly review the basic facts and notation from Sullivan's theory. A \emph{$\Lambda$-algebra} is a commutative differential graded algebra (cdga) of the form $(\land V,d)$, where $V= V^{\geq 0}$ is a graded vector space and $\land V$ is the free graded commutative algebra generated by $V$. Here the differential is required to satisfy the \emph{Sullivan condition}: $V = \cup_{n\geq 0} V(n)$, where
\begin{eqnarray}\label{l1}
V(0) = V\cap \mbox{ker}\, d \hspace{5mm}\mbox{and } V(n+1) = V\cap d^{-1}(\land V(n)).\end{eqnarray}
Following our convention, when there is no ambiguity we suppress the differential from the notation and write
$$\land V = (\land V,d).$$

The \emph{cylinder object} for a cdga $A$ is the cdga $A\otimes \land (t,dt)$ in which deg$\, t= 0$, together with the morphisms $\varepsilon_0, \varepsilon_1: A\otimes \land (t,dt) \to A$ sending $t\mapsto 0,1$. Then two morphisms $\varphi_0, \varphi_1: \land V\to A$ are \emph{homotopic} if some morphism $\varphi : \land V\to A\otimes \land (t,dt)$ satisfies $\varepsilon_i\circ \varphi = \varphi_i$. This is denoted by $\varphi_0\sim \varphi_1$. In particular, a quasi-isomorphism between $\Lambda$-algebras satisfying $H^0= \mathbb Q$ is a homotopy equivalence.

Now $\land V  = \oplus_{p\geq 0} \land^pV$, where $\land^pV$ denotes the linear span of the monomials in $V$ of length $p$; $p$ is called the \emph{wedge degree}. In particular, a $\Lambda$-algebra is \emph{minimal} if $d : V\to \land^{\geq 2} V$ and \emph{quadratic} if $d : V\to \land^2V$. Thus   a minimal $\Lambda$-algebra $(\land V,d)$ determines the associated quadratic $\Lambda$-algebra $(\land V,d_1)$ defined by: $d_1v$ is the component of $dv$ in $\land^2V$.
 
 A \emph{Sullivan algebra} is   a $\Lambda$-algebra $(\land V,d)$ with $V= V^{\geq 1}$. A (minimal) Sullivan model for a cga $A$ is a quasi-isomorphism $\land V \stackrel{\simeq}{\to} A$ from a (minimal) Sullivan algebra, and if $H^0(A)= \mathbb Q$ then $A$ has a minimal Sullivan model. Moreover, given Sullivan models for $A$ and $B$, a morphism $\sigma : A\to B$ determines a unique homotopy class of morphisms, $\psi_\sigma$, for which
 $$
 \xymatrix{
 \land V \ar[d]^{\psi_\sigma}\ar[rr]^\simeq && A \ar[d]^\sigma \\
 \land W \ar[rr]^\simeq && B}$$
 is homotopy commutative; $\psi_\sigma$ is a \emph{Sullivan representative} for $\sigma$. Finally $\Lambda$-algebras can be reduced to minimal Sullivan algebras: if $\land W$ is a $\Lambda$-algebra and $H^0(\land W)= \mathbb Q$ then $\land W \cong (\land V,d) \otimes \land (U\oplus dU)$, where $(\land V,d)$ is a minimal Sullivan algebra and $d : U \stackrel{\cong}{\rightarrow} dU$. Thus the inclusion $\land V\to \land W$ is a homotopy equivalence, as is the surjection $\land W\to \land V$.

 Sullivan algebras link topological spaces and simplicial sets with their   Sullivan completions via two pairs of adjoint functors
 $$Sing : \mbox{Top} \leadsto \mbox{Simp} \hspace{5mm} \mbox{and } \vert\,\,\vert \,:\, \mbox{Simp}\leadsto \mbox{Top},$$
 and
 $$A_{PL}= \mbox{Simp}(-, {A_{PL}}_*) : \mbox{Simp}\leadsto \mbox{Cdga} \hspace{5mm}\mbox{and } \langle\,\,\, \rangle = \mbox{Cdga}(-, {A_{PL}}_*): \mbox{Cdga} \leadsto \mbox{Simp}.$$
 (Here Simp is the category of simplicial sets and $Sing\, X$ is the simplicial set of singular simplices on $X$, and ${A_{PL}}_*$ is the simplicial cdga with $(A_{PL})_n$ the rational cdga generated by the coordinate functions on $\Delta^n$.)
   For simplicity we write   $A_{PL}(X) := A_{PL}(Sing\, X)$ for topological space, $X$.  
 
 \vspace{3mm} For any $\Lambda$-algebra $\land V$, $id_{\langle\land V\rangle}$ is adjoint to a morphism
 $$m_{\land V} : \land V\to A_{PL}\langle \land V\rangle.$$

 Moreover, the functor $\langle \, \rangle$ associates to each morphism of $\Lambda$-algebras,  $\varphi: \land V\to \land W$, a morphism of simplicial sets   $ \langle \land W\rangle \to \langle \land V\rangle,$ and
  $$A_{PL}\langle \varphi\rangle \circ m_{\land V} = m_{\land W}\circ \varphi.$$ 
It is immediate from the definition that $\langle\, \rangle$   has the following properties
\begin{eqnarray}
\label{l2}
\left.
\begin{array}{ll}
(i) & \langle \land V\otimes \land W\rangle = \langle \land V\rangle \times \langle \land W\rangle\\
(ii) & \mbox{If $\varphi_0\sim \varphi_1: \land V\to \land W$ then $\langle \varphi_0\rangle $ and $\langle \varphi_1\rangle$ are homotopic.}\\
(iii) &  \langle\, \rangle  \, \mbox{converts direct limits to inverse limits and}\\ & \mbox{quasi-isomorphisms to homotopy equivalences}.
\end{array}\right\}
\end{eqnarray}
This yields

\vspace{3mm} From the adjoint functors above it follows that, for a simplicial set $X$ and a $\Lambda$-algebra $\land V$,  adjoint to any morphism $\varphi : \land V\to A_{PL}(X)$ is  a simplicial map
$$\widetilde{\varphi} :   X \to \langle \land V\rangle.$$

On the other hand, suppose $\varphi : \land V\to A_{PL}(X)$ is any cdga morphism from a Sullivan algebra. Then $\varphi$ and $m_{\land V}$  are connected via the commutative diagram
\begin{eqnarray}
\label{l3}
\xymatrix{ A_{PL}(X)  && A_{PL}\langle\land V\rangle\ar[ll]^{A_{PL}(\widetilde{\varphi})}\\
 & \land V\ar[lu]^\varphi \ar[ru]_{m_{\land V}}}
 \end{eqnarray}

 \vspace{3mm}\noindent {\bf Definition.} \begin{enumerate}
 \item[(i)]  If $A$ is a cdga then $\langle A\rangle$ is its \emph{simplicial realization}.
 \item[(ii)] A \emph{(minimal) Sullivan model} for a  connected space   $X$, is a quasi-isomorphism from a (minimal) Sullivan algebra,
 $$\xymatrix{\varphi: \land V\ar[rr]^\simeq && A_{PL}(X)}.$$
 \item[(iii)]  For such a minimal Sullivan model,
 $$\xymatrix{\widetilde{\varphi}: X\ar[rr]&& X_{\mathbb Q}:=\langle \land V\rangle}$$
 is the \emph{Sullivan completion} of $X$.
 \item[(iv)] A \emph{Sullivan representative} for a map $f : X\to Y$ is a Sullivan representative for $A_{PL}(f)$. Its simplicial realization, $f_{\mathbb Q}: X_{\mathbb Q}\to Y_{\mathbb Q}$, satisfies $\widetilde{\varphi}_Y\circ f \sim f_{\mathbb Q}\circ \widetilde{\varphi}_X$. 
 \end{enumerate}

\vspace{3mm}\noindent {\bf Example.}
 If each dim$\, V^i<\infty$ then $H(\land V)$ is a graded vector space of finite type, and by (\cite[Theorem 5.4]{RHTII}) $m_{\land V}$ is a quasi-isomorphism. Therefore, in this case
 \begin{eqnarray}
 \label{l4}
 H(\land V)= H(\langle \land V\rangle;\mathbb Q) \hspace{3mm}\mbox{and } \langle \land V\rangle = \langle \land V\rangle_{\mathbb Q}.
 \end{eqnarray}

 \vspace{3mm}
 Next, recall (\cite[Theorem 3.1]{RHTII}) that any morphism $\varphi : \land V\to \land W$ of Sullivan algebras factors as
 $$\land V\to \land V\otimes \land Z\stackrel{\simeq}{\to} \land W,\hspace{3mm} v\mapsto v\otimes 1,$$
 in which $Z = \cup_{k\geq 0} Z(k)$ and $d: 1\otimes Z(k+1) \to \land V\otimes \land Z(k)$. The inclusion $\land V\to \land V\otimes \land Z$ is called a  \emph{$\Lambda$-extension}, and   $\land V\otimes \land Z= \land (V\oplus Z)$  and the quotient $(\land Z, \overline{d}):= \mathbb Q\otimes_{\land V}(\land V\otimes \land Z)$ are $\Lambda$-extensions. By (\cite[Proposition 17.9]{FHTI}) the sequence
 $$\langle\land V\rangle \leftarrow \langle \land V\otimes \land Z\rangle \leftarrow \langle \land Z\rangle$$
 is a Serre fibration and  $\langle \land V\otimes \land Z\rangle$ is homotopy equivalent to $\langle \land W\rangle$.
  
 In particular, when $\varphi$ is a Sullivan representative for $f: X\to Y$ then $\langle \land Z\rangle$ is homotopy equivalent  to the homotopy fibre of $F_{\mathbb Q}: X_{\mathbb Q}\to Y_{\mathbb Q}$. \emph{However it is important to note that $\langle \land Z\rangle$ may not be the Sullivan completion of a connected space.}

Finally, any simplicial map $\sigma: X\to \langle \land V\rangle$ is the adjoint of a morphism $\varphi : \land V\to A_{PL}(X)$. If $X$ is connected and $\psi: \land W\to A_{PL}(X)$ is a  Sullivan model, then $\varphi$ lifts up to homotopy through $\psi$ to give a morphism $\chi : \land V\to \land W$. It is immediate from the definition that the diagram
 \begin{eqnarray}
 \label{l5}
 \xymatrix{X_{\mathbb Q}=\langle \land W\rangle\ar[rr]^{\langle  \chi\rangle} && \langle \land V\rangle\\
 & X\ar[ul]^{\widetilde{\psi}}\ar[ur]_\sigma}
 \end{eqnarray}
 homotopy commutes.

\section{Homotopy groups, examples}

Let $\land V$ be a minimal Sullivan algebra. Then    the adjoint functors in \S 2, together with the fact that $H(S^n)$ and $A_{PL}(S^n)$ have isomorphic minimal Sullivan models, yield natural set bijections (\cite[Theorem 1.3]{RHTII}),
\begin{eqnarray}
\label{l6}
\pi_n\langle \land V\rangle \stackrel{\cong}{\longrightarrow}  (V^n)^\vee\,, \hspace{1cm} n\geq 1.
\end{eqnarray}
These  are isomorphisms of abelian groups when $n\geq 2$ (\cite[Theorem 1.4]{RHTII}), and so  identify $\pi_{\geq 2}\langle \land V\rangle$ as a rational vector space.   The group multiplication in  $\pi_1\langle \land V\rangle$ also has an explicit description in terms of $\land V$   (cf \S 4).
 
 Then, if $\varphi: \land V\to \land W$ is a morphism of minimal Sullivan algebras, a linear map $\varphi_0: V\to W$ is defined by $\varphi_0v$ is the component of $\varphi v$ in $W$. It is immediate from the construction of the bijections 
(\ref{l6}) that they identify 

\begin{eqnarray}
\label{nl7}
\pi_*(\widehat{\varphi}) = \varphi_0^\vee.
\end{eqnarray}
In particular there follows

\begin{Prop}
\label{p1}
Suppose $\varphi : \land V\to \land W$ is a morphism of $\Lambda$-algebras for which $H^0(\land V)= \mathbb Q= H^0(\land W)$. Then the following conditions are equivalent
\begin{enumerate}
\item[(i)] $H(\varphi)$ is an isomorphism
\item[(ii)] $\langle \varphi\rangle$ is a homotopy equivalence.
\end{enumerate}
\end{Prop}

\noindent {\sl proof:} Because these  $\Lambda$-algebras are the tensor product of an acyclic $\Lambda$-algebra and a minimal Sullivan algebra, it is sufficient to consider the case both are minimal. In this case if $H(\varphi)$ is an isomorphism $\varphi$ must be an isomorphism and $\langle \varphi\rangle$ is a homeomorphism.

In the other direction, suppose $\langle \varphi\rangle$ is a homotopy equivalence. Then it follows from (\ref{nl7}) that $\varphi_0$ is an isomorphism and hence $\varphi$ is an isomorphism. \hfill$\square$

\vspace{3mm}\noindent {\bf Example.} Note that $(S^1\vee S^2)_{\mathbb Q}\neq S^1_{\mathbb Q}\vee S^2_{\mathbb Q}$. Indeed let $\land V$ be the minimal Sullivan model of $S^1\vee S^2$. Then dim$\, V^2=\infty$ and so $\pi_2((S^1\vee S^2)_{\mathbb Q})$ is an uncountable vector space. On the other hand, since the universal cover of $S^1_{\mathbb Q}\vee S^2_{\mathbb Q}$ is a wedge of countably many rational spheres $S^2_{\mathbb Q}$, $\pi_2(S^1_{\mathbb Q}\vee S^2_{\mathbb Q})$ is a countable vector space.

\vspace{3mm}\noindent {\bf Example.} For a path connected space $X$ with minimal Sullivan model $(\land V,d)$, by \cite[Theorem 9.2]{RHTII}, cat$(\land V)\leq $ cat$\, X$, where cat$(\land V)$ is the Lusternik-Schnirelmann category of $\land V$  defined in \cite[\S 9.1]{RHTII}. Since the natural map $m_{\land V}: \land V\to A_{PL}\langle \land V\rangle$ has a homotopy section (\cite[Theorem 1.13]{RHTII}), we have also
$$\mbox{cat}(\land V) \leq \mbox{cat}\, \langle \land V\rangle.$$

 \vspace{3mm}  $\Lambda$-algebras, $\land V$, are equipped with an additional structure:   the family,
${\mathcal I}_V = \{V_\alpha\},$
  of finite dimensional subspaces $V_\alpha \subset V$ for which $\land V_\alpha$ is preserved by $d$. For simplicity, we write
  $${\mathcal I}_V = \{\alpha\}.$$
It follows from the Sullivan condition that ${\mathcal I}_V$ is closed under arbitrary intersections and finite    sums. In particular, ${\mathcal I}_V$ is a directed set under inclusion, and 
  \begin{eqnarray}\label{l7} V = \varinjlim_\alpha V_\alpha.\end{eqnarray}  
  Because $\langle\,\,\rangle$ and $(\, )^\vee$ convert  direct limits to inverse limits,  this gives
  \begin{eqnarray}
  \label{l8}
  \langle \land V\rangle = \varprojlim_\alpha \langle \land V_\alpha \rangle \hspace{5mm}\mbox{and } \pi_*\langle \land V\rangle = \varprojlim_\alpha \pi_*\langle \land V_\alpha\rangle.
  \end{eqnarray}

  If $\land V$ is a minimal Sullivan model of $X$, we write $X_\alpha = \langle \land V_\alpha \rangle$ and then
  $$X_\mathbb Q = \varprojlim_\alpha X_\alpha,\hspace{3mm} H(X)= \varinjlim_\alpha H(X_\alpha), \hspace{3mm}\mbox{and } 
  \pi_*(X_{\mathbb Q})= \varprojlim_\alpha \pi_*(X_\alpha).$$
  Moreover, by (\ref{l4}), $
  X_\alpha = (X_\alpha)_{\mathbb Q}.$
  
    \vspace{3mm}\noindent {\bf Example.} A minimal Sullivan algebra $\land V$ has a countable basis if and only if $H(\land V)$ has a countable basis, and this condition holds if and only if in the filtration (\ref{l1}) each $V[n]:= V^{\leq n}(n)$   is finite dimensional. In this case $\{V[n]\}$ is cofinal with ${\mathcal I}_V$.
  In particular, if $\land V$ is the minimal Sullivan model of a space $X$ then this condition holds if and only if each $H^n(X)$ is finite dimensional.

   \vspace{3mm} Next, suppose $\varphi : \land W\to A_{PL}(BG)$ is the minimal Sullivan model of the classifying space of a discrete group $G$. Then $\varphi$ restricts to a morphism 
 $$\varphi_1 : \land W^1\to A_{PL}(BG).$$
 
 \vspace{3mm}\noindent {\bf Definition.} The group homomorphism
 $$\pi_1\langle \varphi_1\rangle : G\to \pi_1\langle \land W^1\rangle$$
 is the \emph{Sullivan completion} of $G$, and we write
 $$G_{\mathbb Q}:= \pi_1\langle \land W^1\rangle.$$
If $\land V$ and $\land W$ are respectively the minimal Sullivan models of $X$ and $B\pi_1\langle X\rangle$ then (\cite[diagram 7.3]{RHTII}) the Sullivan algebras $\land V^1$ and $\land W^1$ coincide. It follows that 
 \begin{eqnarray}\label{l9}
 \pi_1(X_{\mathbb Q}) = [\pi_1(X)]_{\mathbb Q}.
 \end{eqnarray}

 \vspace{3mm}\noindent {\bf Example.} Let $X$ be the wedge of two circles, $X = S^1\vee S^1$. Then $X_{\mathbb Q}=K(G,1)$ where $G$ is the
Sullivan completion of a free group on two generators.  Since $S^1$ is a retract of $X$, $S^1_{\mathbb Q}$ is a retract of $X_{\mathbb Q}$ and by \cite[Example 4.10]{cat}, cat$(X_{\mathbb Q}) \geq 2$.

\vspace{3mm}  A minimal Sullivan algebra also provides an expression for the \emph{Hurewicz homomorphism},
  $$hur : \pi_*\langle \land V\rangle \to H_*(\langle \land V\rangle;\mathbb Q).$$ 
In fact,  division by $\land^{\geq 2}V$ induces a linear map $H^{\geq 1}(\land V) \stackrel{\xi}{\to} V$. Dualizing gives a linear map $\xi^\vee: \pi_*\langle \land V\rangle \to  (H(\land V))^\vee$. If dim$\, V<\infty$ then the inclusion $H_*(\langle \land V\rangle;\mathbb Q)\to    (H^*\langle \land V\rangle)^\vee  $ is an isomorphism. In this case,  by \cite[Proposition 1.19]{RHTII},  $\xi$ is the classical Hurewicz homomorphism.
  
For general minimal Sullivan algebras note that dualizing $H(m_{\land V})$  yields the linear map
  $$j: H_*(\langle \land V\rangle;\mathbb Q) \to H(\langle \land V\rangle)^\vee\to H(\land V)^\vee = \varprojlim_\alpha H(\land V_\alpha)^\vee = \varprojlim_\alpha H_*(\langle \land V_\alpha \rangle;\mathbb Q).$$ 
    Here we have the commutative diagram 
  $$\xymatrix{
  \pi_*\langle \land V\rangle \ar[rr]^\cong \ar[rrd]_{\xi^\vee}\ar[d]^{hur} && \displaystyle\varprojlim_{\alpha \in {\mathcal I}_V} \pi_*\langle \land V_\alpha\rangle\ar[d]^{\varprojlim_\alpha hur_\alpha}\\
  H_*(\langle \land V\rangle;\mathbb Q)\ar[rr]_j && \displaystyle\varprojlim_{\alpha \in {\mathcal I}_V} H_*(\langle \land V_\alpha\rangle;\mathbb Q).}$$

Finally, the homology $H_*(X_{\mathbb Q};\mathbb Z)$ remains mysterious, as it is even unknown   whether or not is is a rational vector space. We do have, as pointed to us by J. Rosenberg, the following.

\begin{Prop} \label{prop2} If in a minimal Sullivan algebra, $\land V$, dim$\,V^1<\infty$,     then $H_{\geq 1}(\langle \land V\rangle;\mathbb Z)$ is a rational vector space.\end{Prop}

\noindent {\sl proof:}  
 Denote $\pi_1\langle \land V\rangle$ by $G$. According to \cite[Proposition 17.9]{FHTI}, $\langle \land V\rangle$ decomposes as a (Serre) fibration
 $$\langle \land V^{\geq 2}\rangle \to \langle \land V\rangle \to \langle \land V^1\rangle.$$
 The Serre spectral sequence for this fibration then converges to $H_*(\langle \land V\rangle;\mathbb Z)$ from the group homology
 $$\mbox{Tor}_p^{\mathbb Z[G]} (\mathbb Z, H_q(\langle \land V^{\geq 2}\rangle;\mathbb Z)).$$
 It is therefore sufficient to show that $\mbox{Tor}_{\geq 1}^{\mathbb Z[G]}(\mathbb Z, \mathbb Z)$ and $H_{\geq 1}(\langle \land V^{\geq 2}\rangle;\mathbb Z)$ are rational vector spaces.
 
 The second assertion is (\cite[Chap II, Proposition 1.1]{Hilton}), since $\pi_*(\langle \land V^{\geq 2}\rangle)$ is a rational vector space. For the first assertion, note that $\langle \land V^1\rangle$ is an Eilenberg-MacLane space $K(G,1)$. Since $G$ is the direct limit of finitely generated nilpotent groups $G_\lambda$, it is sufficient to show that the image of $H_{\geq 1}(K(G_\lambda, 1);\mathbb Z)$ in $H_{\geq 1}(K(G,1);\mathbb Z)$ is contained in an abelian subgroup which is itself a rational vector space.
 
 But since dim$\, V^1<\infty$, $\langle \land V^1\rangle= \langle \land (v_1, \dots, v_r)\rangle$. Since $\langle \land v_i\rangle = S^1_{\mathbb Q}$ and $H_{\geq 1}(S^1_{\mathbb Q}, \mathbb Z)= \mathbb Q$, it follows by induction via a Serre spectral sequence argument that $H_{\geq 1}(\land (v_1, \dots , v_r)\rangle;\mathbb Z)$ is a rational vector space.
 \hfill$\square$

  \vspace{3mm} Finally, even in the simply connected case, the homology of a Sullivan completion can be enormous, as the following Proposition shows.

\begin{Prop} \label{prop3}  If $(\land V,d)$ is a minimal Sullivan algebra in which $H^{\geq 1}(\land V)= H^2(\land V)$ is countably infinite, then  $\,H_3(\langle \land V,d\rangle) $ is uncountable.\end{Prop}

\vspace{3mm}\noindent {\bf Corollary}. If $X$ is a countable wedge of 2-spheres then $H_3(X_{\mathbb Q};\mathbb Q)$ is uncountable.
   
 \vspace{3mm}\noindent {\sl proof of Proposition \ref{prop3}:} $H^2(\land V) $ is the union of an increasing sequence of finite dimensional vector spaces $A(n)$. Let $(\land V(n),d)$ be the minimal model of the cdga  $(\mathbb Q\oplus A(n),0)$. Denote by $L$ the homotopy Lie algebra of $(\land V,d)$ and by $L(n)$ the homotopy Lie algebras for the $(\land V(n),d)$.
 
 Since $A(n)$ is finite dimensional, by \cite[Example 2, \S 18]{FHTI}) $L(n) $ is the free graded Lie algebra $\mathbb L(L(n)_1)$. On the other hand,   $L = \varprojlim_n L(n)$.
 Thus the natural map   $$\varphi : \mathbb L (L_1) \to L$$ is injective and  an isomorphism in degree $1$.

Now decompose a basis of $V^2$   into two infinite sequences,   $x_k, k\geq 1$ and   $y_k, k\geq 1$. There are then elements $z_{ij}$, $t_{ij}$ and $t'_{ij}$ in $V^3$ with
 $$dz_{ij} = x_iy_j\,,\, dt_{ij} = x_ix_j,\,\, dt'_{ij}= y_iy_j.$$
 
 Then let $W$ be the linear span of the $z_{ii}$. The dual elements $f_i\in W^\vee$ defined by $f_i(z_{jj})= \delta_{ij}$ are the basis of a countable subspace of $W^\vee$; fix a direct summand $Q\subset W^\vee$ of this subspace. Then choose elements $u_\alpha\in L_2$ so that the $su_\alpha$ vanish on the $t_{ij}$, the $t'_{ij}$ and on the $z_{ij}$, $i\neq j$, and so that when restricted to $W$ the $su_\alpha$ are a basis of $Q$. Thus the $su_\alpha$ are an uncountable set, and no linear combination of the $su_\alpha$ can vanish on all but finitely many $z_{ii}$. In particular
 the $u_\alpha$ represent linearly independent elements   in $L_2/ \mathbb L^2(L_1)$
 
 \vspace{3mm} Next set $X= \langle \land V\rangle$, and let   $S$ be a wedge of rational spheres $S^2_{\mathbb Q}$ admitting a map $f:S \to X$ that induces an isomorphism on $\pi_2 $. Then $\pi_{\leq 2}(X,S)= 0$, and  
 $$  \pi_3(S)\to \pi_3(X)\to \pi_3(X,S)\to 0$$
 is exact. 
Denote by $\overline{su_\alpha}$ the image of  $su_\alpha$ in $\pi_3(X,S)$. Since the image of $\pi_3(S)$ in $\pi_3(X)$ is contained in $s\mathbb L^2(L_1)$,   the $\overline{su_\alpha}$ are linearly independent.
 
 Next, denote by $C$ the homotopy cofiber of $f$, $C = X \cup_S CS$, where $CS$ is the cone on $S$. Then  $$\pi_r(CS,S)\stackrel{\cong}{\longrightarrow} \pi_{r-1}(S).$$
 
Thus the pairs $(X,S)$ and $(CS,S)$ are both $2$-connected. It follows from the Blakers-Massey theorem (\cite{BM}) that the natural map
 $$\pi_3(X, S)\to \pi_3(C,CS)$$
 is an isomorphism. From the commutative diagram
 $$\xymatrix{
 \pi_3(X) \ar[rr]^g \ar[d] && \pi_3(C)\ar[d]^\cong\\
 \pi_3(X,S) \ar[rr]^\cong && \pi_3(C, CS)}$$
 it follows that   the $g(\overline{su_\alpha})$ are linearly independent.
 
 Now remark that $C$ is $2$-connected, and that the Hurewicz theorem yields a commutative diagram
 $$\xymatrix{ \pi_3(X)\ar[d] \ar[rr]^g && \pi_3(C)\ar[d]^\cong\\
 H_3(X) \ar[rr] &&H_3(C).}$$
 It follows that the images under the Hurewicz map of the elements $su_\alpha$ are linearly independent in $H_3(X;\mathbb Q)$, and so $H_3(X;\mathbb Q)$ is uncountable.\hfill$\square$

  \section{The homotopy Lie algebra, $L_V$, and the completions of $UL_V$ and $H_*(\Omega X_{\mathbb Q};\mathbb Q)$}
  
 For each minimal Sullivan algebra $\land V$, the associated quadratic differential, $d_1$,  determines the \emph{homotopy Lie algebra}, $$L_V :=\{(L_V)_p\}_{p\geq 0}.$$ $L_V$ is defined by 
  a (degree 1) suspension isomorphism $s:(L_V)_p \stackrel{\cong}{\rightarrow}  (V^{p+1})^\vee$ and
  $$\langle v, s[x,y]\rangle = (-1)^{deg\, y+1} \langle d_1v,sx,sy\rangle\,, \hspace{1cm} v\in V.$$

  Moreover, the family ${\mathcal I}_V = \{V_\alpha\}$   endows $L_V$ with additional structure. Denote by $L_\alpha$ the homotopy Lie algebra of $\land V_\alpha$. Since dim$\, V_\alpha<\infty$, $L_\alpha$ is a nilpotent graded Lie algebra.
  The duals of the inclusions $V_\alpha \to V$ then desuspend to Lie algebra surjections   $\rho_\alpha :L_V\to L_\alpha$, which define an isomorphism
  $$L_V\stackrel{\cong}{\longrightarrow} \varprojlim_\alpha L_\alpha.$$
  
On the other hand, recall that the classical completion of the universal enveloping algebra, $UL$, of a graded Lie algebra $L$ is the inverse limit
$$\widehat{UL} = \varprojlim_n \, UL/I^n,$$
$I^n$ denoting the $n^{th}$ power of the augmentation ideal $I$ generated by $L$. It turns out that if $L_V$ is the homotopy Lie algebra of a minimal Sullivan algebra, $\land V$,  then a more useful notion is the \emph{Sullivan completion}, $\overline{UL_V}$, given by
$$\overline{UL_V}= \varprojlim_{\alpha \in {\mathcal I}_V} \, \widehat{UL_\alpha}.$$
In general, the natural map $\widehat{UL_V}\to \overline{UL_V}$ may not be an isomorphism, but it is an isomorphism   when $L_V$ is finitely generated.

  \vspace{3mm}\noindent {\bf Remark}. Let $x_i$ be a graded basis for $UL_\alpha$; then $\widehat{UL_\alpha} = \prod_i \mathbb Q \cdot x_i$. It follows that $\overline{UL_V} = \prod_j \mathbb Q \cdot y_j$, where $y_j$ ia a graded basis for $UL$.
  
  \vspace{3mm}
  
  The isomorphism $L_V = \varprojlim_\alpha L_\alpha$ also identifies the product structure in $\pi_1\langle \land V\rangle$. 
Since dim$\, L_\alpha<\infty$, bijections
$$\exp_\alpha : (L_\alpha)_0 \to G_\alpha \subset \widehat{UL_\alpha}$$
onto the group $G_\alpha$ of units of $\widehat{UL_\alpha}$ are given by $x\mapsto \sum_{n\geq 0} x^n/n!$ (\cite[Chapter 4]{Serre}, \cite[Chapter 2]{RHTII}). The inverse bijection $\log_\alpha : G_\alpha \to  (V_\alpha^1)^\vee$ is then given by the standard power series. Moreover, by  (\cite[Theorem 2.4]{RHTII})  the composites
$$\pi_1\langle \land V_\alpha\rangle \to (V_\alpha^1)^\vee\to (L_\alpha)_0 \to G_\alpha$$
are isomorphisms $\pi_1\langle V_\alpha \rangle \stackrel{\cong}{\longrightarrow} G_\alpha$ of groups. Thus passing to inverse limits yields a group isomorphism,
$$\exp: \pi_1\langle \land V\rangle \stackrel{\cong}{\longrightarrow} \varprojlim_{\alpha \in {\mathcal I}_V}\, G_\alpha \subset \overline{UL_V},$$
whose inverse bijection is $\log = \varprojlim_\alpha \log_\alpha$.

In particular, for $x\in (L_\alpha)_0$, $\exp (px) = (\exp x)^p$, $p\in \mathbb N$, it follows that $a\mapsto a^p$ is a bijection in each $G_\alpha$. Therefore
$$x\mapsto x^p$$
is a bijection in $\pi_1\langle \land V\rangle$; i.e., $\pi_1\langle \land V\rangle$ is \emph{uniquely divisible}.

 Finally, for any minimal Sullivan algebra $\land V$ the Whitehead product
 $$[\,,\,]_W: \pi_p\langle \land V\rangle \, \times \, \pi_q\langle \land V\rangle \to \pi_{p+q-1}\langle \land V\rangle$$
 may be computed directly from $L_V$ as follows  
\begin{eqnarray}
\label{l10}
[sx,sy]_W= \left\{
 \renewcommand{\arraystretch}{1.3}
 \begin{array}{ll}
 (-1)^{deg\, x} s[x,y], \hspace{2cm} & x,y\in (L_V)_{\geq 1}\\
 s\, \mbox{Ad}(\exp x)y-sy, & x\in (L_V)_0, y\in (L_V)_0\\
 s(\log  [\exp x,\exp y]), & x,y\in (L_V)_0.
 \end{array}
 \renewcommand{\arraystretch}{1}
 \right.\end{eqnarray}
This follows from \cite[Chapter 2 and Theorem 4.2]{RHTII} for each $\land V_\alpha$, and then by an inverse limit argument.

\vspace{3mm}\noindent {\bf Remark:} The family of finite dimensional subspaces $W\subset V$ for which $\land W$ is preserved by $d_1$ is cofinal with ${\mathcal I}_V$. It follows that $\widehat{UL_V}$, and the Whitehead products including the group structure in $\pi_1\langle \land V\rangle$ depend only on the associated quadratic Sullivan algebra $(\land V,d_1)$.

\vspace{3mm}\noindent {\bf Example.} Let $(\land V,d)$ be the  minimal Sullivan model of $X= S^1\vee S^3$. Then $V^1= \mathbb Q v$ and $V^3$ has a basis $y_0, 
\dots , y_n, \dots$
satisfying $dy_0=0$ and for $n\geq 1$, $dy_n = y_{n-1}v$. Choose $x\in \pi_1(X)$ with $\langle v,x\rangle = 1$. We identify $\pi_3\langle \land V, d\rangle =  (V^3)^\vee$ with the space of series $\mathbb Q[[t]]$ by associating to a linear map $f$ the series $g(t)=\sum \langle y_n, f\rangle t^n$. Using this identification, by \cite[Theorem 4.6]{RHTII}, the action of $\pi_1\langle \land V\rangle$ on $\pi_3\langle \land V\rangle$ is given as the product of two series: 
$$(\alpha \bullet g)(t) = \exp (t\log \alpha) \cdot g(t)\, \hspace{1cm} \alpha \in \pi_1\langle \land V^1\rangle.$$

\vspace{3mm}
Next, in any minimal Sullivan algebra, $\land V$, each $\alpha \in {\mathcal I}_V$ yields the map
$$f_\alpha :   \langle \land V\rangle \to \langle \land V_\alpha\rangle,$$
and if $\alpha \leq \beta$ then $f_\alpha = f_{\alpha \beta}\circ f_\beta$ where $f_{\alpha\beta}: \langle \land V_\beta\rangle \to \langle \land V_\alpha \rangle$ is induced by the inclusion. This defines     a map
$$H_*(\Omega \langle \land V\rangle;\mathbb Q)\to   \varprojlim_\alpha H_*(\Omega\langle \land V_\alpha\rangle;\mathbb Q).$$

On the other hand, completing $H_*(\Omega\langle \land V_\alpha\rangle;\mathbb Q)$ with respect to the augmentation ideal yields natural isomorphisms (\cite[Proposition 3.3]{De})
$$\widehat{H}(\Omega\langle \land V_\alpha\rangle;\mathbb Q) \cong \widehat{UL_\alpha}.$$
Thus we obtain the morphism
$$H_*(\Omega \langle \land V\rangle;\mathbb Q)   \to \overline{UL_V} := \varprojlim_\alpha \widehat{UL_\alpha},$$
which identifies $\overline{UL_V}$ as a completion  of $H_*(\Omega \langle \land V\rangle;\mathbb Q)$. We denote this by
$$\overline{H_*(\Omega \langle \land V\rangle;\mathbb Q)} := \overline{UL_V}.$$

Moreover, the isomorphisms $s: L_V\stackrel{\cong}{\longrightarrow} \pi_*\langle \land V\rangle$ and $s: \pi_*\Omega \langle \land V\rangle \stackrel{\cong}{\longrightarrow} \pi_*\langle \land V\rangle$ define an isomorphism
$$L_V \stackrel{\cong}{\to} \pi_*\Omega \langle \land V\rangle.$$
In view of (11), this converts the Lie bracket in $L_V$ to the Samelson bracket in $\pi_*\Omega\langle \land V\rangle $, in degrees $\geq 1$, and thus extends the Samelson bracket to all of $\pi_*(\Omega \langle \land V\rangle).$

    \vspace{3mm}\noindent {\bf Remark.} If $\land V$ is the minimal Sullivan model of a path connected space $X$, then a homological invariant of $X_{\mathbb Q}=\langle \land V\rangle$ provides a lower bound for the   Lusternik-Schnirelmann category, cat$\, X$, of $X$ (\cite{De}),
$$\mbox{Ext}^p_{UL_V}(\mathbb Q, \overline{H_*(\Omega X_{\mathbb Q})})\neq 0, \hspace{5mm} \mbox{some $p\leq \mbox{cat}\, X$}.$$

  \section{The Sullivan topology}
  
  Fix a minimal Sullivan algebra, $\land V$, and recall that ${\mathcal I}_V = \{\alpha\}$ is the index set for the finite dimensional subspaces $ V_\alpha \subset V$ for which $\land V_\alpha$ is preserved by $d$. Then by (\ref{l6}) we may identify
  $$\pi_*\langle \land V,d\rangle = V^\vee \hspace{3mm} \mbox{and } \pi_*\langle \land V_\alpha\rangle= V_\alpha^\vee.$$
  Thus the inclusions $V_\alpha \to V$ dualize to surjections
  $$\rho_\alpha : \pi_*\langle \land V\rangle \to \pi_*\langle \land V_\alpha\rangle.$$

  Moreover, since ${\mathcal I}_V$ is closed under arbitrary intersections and finite sums, the set of subspaces $\{ \mbox{ker}\, \rho_\alpha\}$ is   closed under   finite intersections and arbitrary sums.
  Thus (similar to an observation of Lefschetz (\cite{Le})) it follows that sets of the form
  $${\mathcal{O}} = \cup_{x_\alpha \in L_V, \alpha \in {\mathcal I}_V }\,\,\, \mbox{ker}\, \rho_\alpha + x_\alpha$$
  are the open sets of a topology in $\pi_*\langle \land V\rangle$. Desuspension then transfers this to a topology in $L_V$.

\vspace{3mm}\noindent {\bf Definition.} The topology just defined is the \emph{Sullivan topology} in $L_V$.

\vspace{3mm}\noindent {\bf Remark.}  The Sullivan topology may be
identified with that introduced by Lefschetz. In fact every finite
dimensional subspace $S\subset V$ induces a surjection $\rho_S :
V^\vee\to S^\vee$. Moreover, each $S$ is a subspace of some $V_\alpha$,
so that $\{\mbox{ker}\, \rho_\alpha\}$ is cofinal with $\{\mbox{ker}\,
\rho_S\}$. But the topology in $V^\vee$ determined by $\{\mbox{ker}\,
\rho_S\}$ is the topology introduced by Lefschetz.

\vspace{3mm} The next Proposition in particular exhibits $(\land V, d_1)$ as the continuous dual of the classical \emph{Cartan-Chevalley-Eilenberg differential graded coalgebra} $(\land sL_V, \partial)$,  defined by $\partial (sx\land sy) = (-1)^{1+ deg\, x} s[x,y]$. 

\begin{Prop}\label{prop4} Let $(\land V,d)$ be a minimal Sullivan algebra. With the notation above:
\begin{enumerate}
\item[(i)] The Whitehead products in $\pi_*\langle \land V\rangle$ are continuous in the Sullivan topology. In particular, $\pi_1\langle \land V\rangle$ is a topological group.
\item[(ii)] The Sullivan topology in $L_V$ coincides with that determined by the associated quadratic Sullivan algebra $(\land V,d_1)$.
\item[(iii)] The Lie bracket $[\,,\,] :  L_V\times  L_V\to  L_V$ is continuous in the Sullivan topology.
  \item[(iv)] $\land^pV$ may be identified with the continuous alternating $p$-linear maps
  $$sL_V\times \dots \times sL_V\to \mathbb Q,$$
  where $\mathbb Q$ has the discrete topology, and the differential, $d$, is continuous.
  \end{enumerate}\end{Prop}
  
 \noindent {\sl proof:} (i) The product in $\pi_1\langle \land V^1\rangle = (V^1)^\vee$ is given
by $x\cdot y = \log (\exp x\cdot \exp y)$, where $\exp x$, $\exp y$ are
elements in the group  of units in $\overline{UL_V}$. Since
$U\rho_\alpha: \overline{UL_V}\to \widehat{UL_\alpha}$ commutes with
$\exp$ and $\log$ it follows that if $\widehat{x}\in \mbox{ker}\,
\rho_\alpha + x$ and $\widehat{y}\in \mbox{ker}\, \rho_\alpha + y$ then
$\widehat{x}\cdot \widehat{y}\in \mbox{ker}\, \rho_\alpha ,+ x\cdot y$.
Thus multiplication is continuous. Similarly inversion is also
continuous. The same argument shows that $\pi_1\langle \land V^1\rangle$
acts continuously in $L_V$. Thus (i) follows from (\ref{l10}).
  
  (ii) This follows from the fact that if $\rho'_\beta: L_V\to L_\beta$ are the surjections determined by $(\land V, d_1)$ then $\{\mbox{ker}\, \rho'_\beta\}$ is cofinal with $\{\mbox{ker}\, \rho_\alpha\}$. 
  
  (iii) This follows because if $[x,y]\in \mbox{ker}\, \rho_\alpha + z$ then $[x+\mbox{ker}\, \rho_\alpha, y+\mbox{ker}\, \rho_\alpha]\subset \mbox{ker}\, \rho_\alpha+ z$.
  
  (iv) If $\Phi\in \land^pV$ then $\Phi\in \land^pV_\gamma$, some $\gamma$. Now $\langle V_\gamma, s\mbox{ker}\, \rho_\gamma\rangle = 0$. Thus, if $x_1, \dots , x_p\in L_V$, then $\langle \Phi, s(x_1+\mbox{ker}\, \rho_\gamma), \dots , s(x_p+\mbox{ker}\, \rho_\gamma)\rangle = \langle \Phi, sx_1, \dots , sx_p\rangle$. This shows that $\Phi$ acts continuously.
  
  On the other hand, suppose $\Phi : sL_V\times \dots \times sL_V\to \mathbb Q$ is any continuous alternating multilinear function. Then for some $\gamma$, $\Phi$ vanishes on $s\,\mbox{ker}\, \rho_\gamma \times \dots \times s\,\mbox{ker}\, \rho_\gamma$. But this implies that $\Phi$ is in the image of $\land^pV_\gamma$ and hence in $\land^pV$. 
  \hfill$\square$
 
\vspace{3mm}\noindent {\bf Corollary.}  Suppose a graded Lie algebra,
$E$, is the homotopy Lie algebra of two quadratic Sullivan  algebras
$(\land V, d_1)$ and $(\land W, d_1)$. If the induced Sullivan
topologies in $E$ coincide, then $(\land V,d_1)\cong (\land W, d_1)$.

\vspace{3mm} However, the following remains an open problem:

\vspace{3mm}\noindent {\bf Problem.} Can a graded Lie algebra $E$ be the
homotopy Lie algebra of two non-isomorphic quadratic Sullivan algebras ?

\vspace{3mm} Note that with additional data, it is possible to construct
a quadratic Sullivan algebra from a graded Lie algebra $E$, and with a
certain finiteness hypothesis, Proposition \ref{prop5} below then shows that this Sullivan
algebra is unique.

First suppose a graded Lie algebra $E= E_{\geq 0}$ is equipped with a family of surjections $\rho_\alpha : E\to E_\alpha$ onto nilpotent and finite dimensional graded Lie algebras. Then the dual of the Cartan-Chevalley-Eilenberg constructions are quadratic Sullivan algebras $\land sE_\alpha^\vee$. If also $E= \varprojlim_\alpha E_\alpha$ and $\{\mbox{ker}\, \rho_\alpha\}$ is closed under finite intersections and arbitrary sums then $\{sE_\alpha^\vee\}$ is a directed set and 
$$\land W:= \varinjlim_\alpha \land sE_\alpha^\vee$$
is a quadratic Sullivan algebra with homotopy Lie algebra $E$. Moreover, $\{sE_\alpha^\vee\}$ is cofinal with ${\mathcal I}_W$, so that the topology in $E$ determined by $\{\mbox{ker}\, \rho_\alpha\}$ is the Sullivan topology.

\vspace{3mm} Finally, recall that the lower central series for a graded Lie algebra $E= E_{\geq 0}$ is the descending sequence of ideals $E^n$ in which $E^n$ is the linear span of iterated commutators $[x_1, [x_2, \dots , x_n]\dots ]]$ of length $n$. 
By definition, $E$ is \emph{pronilpotent} if $E\stackrel{\cong}{\to} \varprojlim_n E/E^n$. 

\begin{Prop}\label{prop5} Suppose $E$ is a pronilpotent Lie algebra and $E/E^2$ is a graded vector space of finite type. Then $E$ is the homotopy Lie algebra of a unique quadratic Sullivan algebra, $\land W$.\end{Prop}

 \noindent {\sl proof:} Since the Lie bracket defines linear surjections $E/E^2\otimes E^n/E^{n+1}\to E^{m+1}/E^{m+2}$ it follows that each $E/E^n$ is a graded vector space of finite type. Since $E/E^n$ is also nilpotent, the dual of the Cartan-Chevalley-Eilenberg construction is a quadractic Sullivan algebra $\land s(E/E_n)^\vee$. It follows that $\land W = \varinjlim_n \land s(E/E^n)^\vee$ is also a quadratic Sullivan algebra and, since $E$ is pronilpotent, $W^\vee = \varprojlim_n s\, E/E^n = sE$. It is immediate that the corresponding Lie bracket in $E$ is the original Lie bracket.

On the other hand, suppose $E$ is also the homotopy Lie algebra of another quadratic Sullivan algebra, $\land V$. Because $E/E^n$ has finite type, there is a subspace $V(n)\subset V$ such that $\langle V(n), sE^n\rangle= 0$ and the induced pairing $V(n) \times s\, E/E^n \to \mathbb Q$ is non degenerate.

These conditions imply that $\land V(n)$ is a sub quadratic Sullivan algebra with homotopy Lie algebra $E/E^n$. In particular, $\land V(n)$ is the dual of the Cartan-Chevalley-Eilenberg complex. But it is immediate that $\land V = \varinjlim_n \land V(n)$, and so $\land V \cong \land W$.
\hfill$\square$

 \section{Lower central series}
 
 \emph{In this section we fix a minimal Sullivan algebra, $\land V$, with associated quadratic differential, $d_1$.}

 Here, and subsequently, we shall rely on the following property   \cite{Bour}: 

\vspace{3mm}\begin{quote}
{\sl Suppose
 $0\to A_\alpha \to B_\alpha \to C_\alpha\to 0$ are exact sequences of morphisms of inverse systems of vector spaces. Then if each dim$\, C_\alpha <\infty$,
 \begin{eqnarray}\label{l11}0 \to \varprojlim_\alpha A_\alpha \to \varprojlim_\alpha B_\alpha \to \varprojlim_\alpha C_\alpha\to 0\end{eqnarray}
  is also exact.}
  \end{quote}

The classical lower central series of a graded Lie algebra $L$ is the     sequence of ideals $L= L^1\supset L^2\supset \dots \supset L^r\supset \dots$ in which $L^r$ is the 
linear span of iterated commutators of length $r$, $[x_1,[\dots [x_{r-1},x_r]\dots ]$. Again, a more natural role in Sullivan's theory is played by the ideals
$$L_V^{(r)} := \varprojlim_{\alpha \in {\mathcal I}_V}\, L_\alpha^r,$$
when $L_V$ is the homotopy Lie algebra of   $\land V$. The sequence 
$$L_V = L_V^{(1)}\supset \dots \supset L_V^{(r)} \supset \dots$$
is the \emph{Sullivan lower central series}. It satisfies
$$L_V^r \subset L_V^{(r)}\,, \hspace{5mm} [L_V^{(r)}, L_V^{(s)}]\subset L_V^{(r+s)}\,, \hspace{3mm}\mbox{and } \cap_r L_V^{(r)}= 0.$$
Further, since each $L_\alpha$ is finite dimensional,
$$L_V^{(r)}/L_V^{(s)} = \varprojlim_{\alpha \in {\mathcal I}_V} \, L_\alpha^r/L_\alpha^s,\hspace{3mm}\mbox{and}$$
\begin{eqnarray}
\label{l12}
L_V =  \varprojlim_\alpha L_\alpha = \varprojlim_\alpha \varprojlim_r L_\alpha /L_\alpha^r = \varprojlim_r \varprojlim_\alpha L_\alpha /L_\alpha^r =\varprojlim_r L_V/L_V^{(r)} .\end{eqnarray}
This exhibits $L_V$ as   complete with respect to the Sullivan lower central series. By contrast the map $L_V\to \varprojlim_r L_V/L_V^r$ may not always be an isomorphism, although we have the

\vspace{3mm}\noindent {\bf Example.} If $H^1(\land V)$ and each $V^n$, $n\geq 2$, are finite dimensional then
$$L_V^{(r)}= L_V^r\,, \hspace{1cm} r\geq 1.$$

\vspace{3mm}Moreover if $\Phi\in \land^{\geq 1}V$ is a $d$-cycle then the component $\Phi_1$ of $\Phi$ in $V$ is a $d_1$-cycle. It follows that the Hurewicz map, $\xi^\vee :sL_V\stackrel{}{\longrightarrow}  H(\land V)^\vee$, described in \S 3 vanishes on $sL_V^2$. Hence $sL^2_{V_\alpha}$ vanishes on $ H(\land V_\alpha)^\vee$. Since $ H(\land V)^\vee= \varprojlim_\alpha  H(\land V_\alpha)^\vee$ it follows that $sL_V^{(2)}$ vanishes on $ H(\land V)^\vee$, so that the Hurewicz map factors to give
\begin{eqnarray}
\label{l13}
s(L_V/L_V^{(2)})\to  H(\land V)^\vee.
\end{eqnarray}

\vspace{3mm} Analogously, the lower central series for a group $G$ is the sequence of normal subgroups $G = G^1\supset \dots \supset G^r\supset \dots $, where $G^{r+1}$ is the subgroup generated by the elements $aba^{-1}b^{-1}$, $a\in G$, $b\in G^r$. Analogous to the construction above, when   $G = \pi_1\langle \land V\rangle$, we set
$$G^{(r)} := \varprojlim_{\alpha \in {\mathcal I}_V} G_\alpha^r,$$
where $G_\alpha = \pi_1\langle \land V_\alpha
 \rangle$. This defines the \emph{Sullivan lower central series}, 
 $$G = G^{(1)}\supset \dots \supset G^{(r)}\supset \dots$$ for $G$. As in the Lie algebra case, it is immediate that
 $$G^r \subset G^{(r)}\,, \hspace{5mm} [G^{(r)}, G^{(s)}]\subset G^{(r+s)}, \hspace{3mm}\mbox{and } \cap_r G^{(r)} = \{e\}.$$

 Moreover (\cite[Theorem 2.4 and Corollary 2.4]{RHTII}), the exponential map restricts to bijections $(L_\alpha)^r_0 \stackrel{\cong}{\longrightarrow} G_\alpha^r$, which then factor to yield linear isomorphisms $(L_\alpha)_0^r /(L_\alpha)_0^{r+1} \stackrel{\cong}{\longrightarrow} G_\alpha^r/G_\alpha^{r+1}$. 
 In particular, these are finite dimensional vector spaces. Passing to inverse limits yields linear isomorphisms
 $$(L_V)_0^{(r)} \, /\, (L_V)_0^{(r+1)} \, \stackrel{\cong}{\longrightarrow}\, G^{(r)}/G^{(r+1)}\,, \hspace{5mm} r\geq 1.$$
 It follows from this, and induction on $s\geq r$ that
 $$G^{(r)}/G^{(s)} \stackrel{\cong}{\longrightarrow} \varprojlim_{\alpha \in {\mathcal I}_V} G_\alpha^r/G_\alpha^s $$
 is an isomorphism of groups, and that when $s=r+1$ these are isomorphisms of rational vector spaces.
 
 Finally, filtering by the normal subgroups $G^{(r)}$ produces an associated graded group
 $$gr\, G = \oplus_{r} G^{(r)}/G^{(r+1)}.$$
 The properties above for the filtration imply that $G$ satisfies  condition N of Lazard \cite{La}, and hence the commutator $[a,b]= aba^{-1}b^{-1}$ makes $gr\, G$ into a graded Lie algebra. On the other hand, the filtration $L_V^{(r)}$ makes $(L_V)_0$ into a filtered Lie algebra, and it is straightforward from the properties above that the associated graded Lie algebra $gr(L_V)_0$ is the Lazard Lie algebra, $gr\, G$:
 $$gr\, G = gr(L_V)_0.$$
 
 \vspace{3mm}\noindent {\bf Example.} \emph{The space $H_1(\langle \land V\rangle)$.}
 
 For each $\alpha \in {\mathcal I}_V$ the isomorphism $H_1(\langle \land V_\alpha\rangle;\mathbb Z) = G_\alpha/ G_\alpha^2 \cong L_\alpha /[L_\alpha, L_\alpha]$ identifies $H_1(\langle \land V_\alpha\rangle; \mathbb Z)$ as the rational vector space $L_\alpha /[L_\alpha, L_\alpha]$. Passing to inverse limits yields the diagram
 $$\xymatrix{
 G/G^2 \ar[rr]^\cong \ar@{->>}[d] && H_1(\langle \land V\rangle;\mathbb Z)\ar[d]\\
 G/G^{(2)} \ar[rr]^\cong && \varprojlim_\alpha H_1(\langle \land V_\alpha\rangle;\mathbb Z).}$$
 Thus if $G^2 = G^{(2)}$ (equivalently if $L^2= L^{(2)}$) then this exhibits $H_1(\langle \land V\rangle;\mathbb Z)$ as a rational vector space.
 
 \vspace{3mm} Now define the \emph{Sullivan completion} of the group ring $\mathbb Q[G]$ by 
 $$\overline{\mathbb Q[G]} := \varprojlim_\alpha \widehat{\mathbb Q[G_\alpha]}$$
 where $\widehat{\mathbb Q[G_\alpha]}$ is the completion of ${\mathbb Q[G_\alpha]}$ with respect to its augmentation ideal. According to  (\cite{Pi}, \cite[Proposition 3.2]{De}), there is  a natural isomorphism
 $$\widehat{\mathbb Q[G_\alpha]} \cong (\widehat{UL_\alpha})_0.$$
 Taking inverse limits gives the isomorphism,
 $$\overline{\mathbb Q[G]}= \varprojlim_\alpha\widehat{\mathbb Q[G_\alpha]} \cong \varprojlim_\alpha (\widehat{UL_\alpha})_0= (\overline{UL})_0.$$
 
\vspace{3mm}  Next observe that  filtration (1) of $V$   determined by $d_1$ and denoted $\{V_n\}$   is given by
 $$V_0= V\cap \mbox{ker}\, d_1 \hspace{5mm}\mbox{and } V_{n+1} = V\cap d_1^{-1}(\land^2V_n), \hspace{3mm} n\geq 0.$$
 Restricting $sL_V$ to $V_n$ gives linear maps $sL_V \to   V_n^\vee$.
 
\begin{Prop}\label{prop6} The linear maps $sL_V\to V_n^\vee$ factor over the surjections $sL_V \to s\left( L_V/L_V^{(n+2)}\right)$ to give degree 1 isomorphisms     
$$
 \label{4}
 L_V\, /\, L_V^{(n+2)} \,\stackrel{\cong}{\longrightarrow}\,  V_n^\vee, \hspace{5mm} n\geq 0,
$$\end{Prop}

\noindent {\bf Corollary.} If $H(\land V, d_1)$ has finite type then dim$\, L_V/L_V^{(n+1)}<\infty$, $n\geq 0$.

\vspace{3mm}\noindent {\sl proof of Proposition \ref{prop6}:} We establish the equivalent assertions,
 \begin{eqnarray}
 \label{l14}
 L_V^{(n+2)} \stackrel{\cong}{\longrightarrow}  (\, V/V(n)\, )^\vee, \hspace{3mm} n\geq 0.
 \end{eqnarray}
 For this, let ${\mathcal I}_{V,1}$ be the set of finite dimensional subspaces $V_\alpha \subset V$ such that $\land V_\alpha $ is preserved by $d_1$. Then let ${\mathcal J}_V$ be the set of subspaces $W\subset V$ for which $\land W$ is preserved by $d$ and dim$\, W\cap\mbox{ker}\, d<\infty$. Note that for any two subspaces $S(1)$ and $S(2)$ of a graded vector space $S$, 
 that $\land S(1)\cap \land S(2)= \land (S(1)\cap S(2))$. In particular, if $V_\alpha\in {\mathcal I}_{V,1}$ and $W
 \in {\mathcal J}_V$ then it follows by induction that $(V_\alpha)_n\cap W = (V_\alpha \cap W)_n$. Hence
 $$L_V^{(n+2)} = \varprojlim_{W\in {\mathcal J}_V} \hspace{2mm} \varprojlim_{W_\alpha \in {\mathcal I}_{W,1}} \hspace{2mm} L_{W_\alpha}^{n+2} = \varprojlim_{W\in {\mathcal J}_V} \hspace{2mm} L_W^{(n+2)}.$$
 
 On the other hand, for $W\in {\mathcal J}_V$ the subspaces $W_k$ are finite dimensional and cofinal in ${\mathcal I}_{W,1}$. Therefore
 $$L_W^{(n+2)} = \varprojlim_k L_{W_k}^{n+2}.$$
 Since $(W(k)_n = W_n$ for $k\geq n$ it follows in this case, exactly as in the proof of Theorem 2.1, p.50 in \cite{RHTII}, that $L^{n+2}_{W_k} \stackrel{\cong}{\longrightarrow}   (W_k\,/\, W_n)^\vee$. This gives
 $$L_W^{(n+2)} = \varprojlim_k   (W_k\,/\, W_n)^\vee =   (W\, /\, W_n)^\vee =   (W\, /\, W\cap V_n)^\vee.$$
 Formula (\ref{l14}) follows. \hfill$\square$
 
 \vspace{3mm} The Hurewicz map (\ref{l13}) for $(\land V, d_1)$ may be regarded as a linear map of degree 1, $\xi^\vee: L_V\to H(\land V, d_1)^\vee$.
 
 \vspace{3mm}\noindent {\bf Corollary.}  The Hurewicz map for $(\land V, d_1)$ factors as
 $$L_V\to L_V/L_V^{(2)} \stackrel{\cong}{\longrightarrow}  V_0^\vee \subset H^{\geq 1}(\land V,d_1)^\vee.$$
 In particular, in degree 1 it translates to
 $$G\to G/G^{(2)} \stackrel{\cong}{\longrightarrow} (H^1(\land V)^\vee,$$
 where $G= \pi_1\langle \land V\rangle$
 
 \vspace{3mm}\noindent {\sl proof:} Because $\xi$ is induced by the surjection $\land^{\geq 1}V\to V$ with kernel $\land^{\geq 2}V$, and because $(\land V, d_1)$ is quadratic, it follows that $\xi$ factors as $\xymatrix{H^{\geq 1}(\land V, d_1)\ar@{->>}[r] & V_0 \ar@{^{(}->}[r] & V}$. This gives the first assertion. The second follows because $H^1(\land V,d)= H^1(\land V, d_1)$. \hfill$\square$
 
 \vspace{3mm}\noindent {\bf Remark.} A second application of the proof of Theorem 2.1 in \cite{RHTII} also gives
 $$L_V^{(n)} = L_V^n\, \hspace{1cm} n\geq 1,$$
 if dim$\, V\cap\mbox{ker}\, d_1<\infty$.

\section{The holonomy representation of a $\Lambda$-extension}

 Now   recall from \S 2 that  
a \emph{$\Lambda$-extension} of a Sullivan algebra $\land V$ is  sequence of cdga morphisms of the form
$$
 \xymatrix{\land V\ar[r]^\lambda & \land V\otimes \land Z \ar[r]^\rho & \land Z}, \hspace{5mm} \lambda v= v\otimes 1,
$$
 in which $(\land Z, \overline{d}) = \mathbb Q\otimes_{\land V}(\land V\otimes \land Z).$

 Now suppose  that $\land V$ is a minimal Sullivan algebra. Then the differential in $\land V\otimes \land Z$ satisfies
\begin{eqnarray}
\label{l15}
d(1\otimes \Phi)= 1\otimes \overline{d} \Phi + \sum_i v_i\otimes \theta_i\Phi + \Omega,\end{eqnarray}
 where $\overline{d}$ is the differential in $\land Z$, $v_i$ is a basis of $V$ and $\Omega \in \land^{\geq 2}V\otimes \land Z$.
 Here the $\theta_i$ are derivations in $(\land Z, \overline{d})$, and setting
 $$\overline{\theta}(x)= -\langle v_i, sx\rangle H(\theta_i), \hspace{1cm} x\in L_V$$
 defines the \emph{holonomy representation} of $L_V$ in $H(\land Z)$.
 
 \vspace{3mm} Next, suppose $ \land V\otimes \land Z$ is a $\Lambda$-extension of a minimal Sullivan algebra,   and set $W = V\oplus Z$. Then $W  $ is the union of the finite dimensional subspaces 
$W_\alpha = V_\alpha \oplus Z_\alpha$ for which $\land V_\alpha\otimes \land Z_\alpha $ is preserved by $d$. 
Under inclusion, these form a directed set ${\mathcal J}$. 
Moreover $\{ V_\alpha, \alpha \in {\mathcal J} \}$, $\{ W_\alpha , \alpha\in {\mathcal J} \}$, and 
$\{ Z_\alpha, \alpha \in {\mathcal J} \})$ are respectively cofinal with ${\mathcal J}_V$, ${\mathcal J}_W$, and ${\mathcal J}_Z$.

In particular, if $\beta\geq \alpha \in {\mathcal J}$ then the $\Lambda$-extension $\land V_\beta\otimes\land Z_\alpha$ yields a holonomy representation of $L_{V_\beta}$ in $H(\land Z_\alpha)$; The defining condition for $\Lambda$-extensions  implies that this representation extends to a representation of $\widehat{UL_\beta}$. Passing to inverse limits then gives a representation of $\overline{UL_V}$ in $H(\land Z_\alpha)$, and passing to direct limits gives a representation of $\overline{UL_V}$ in $H(\land Z)$. 

\vspace{3mm}\noindent {\bf Definition.} This extension of the holonomy representation of $L_V$ to $\overline{UL_V}$ is \emph{the holonomy representation} of $\overline{UL_V}$ in $H(\land Z)$.

\vspace{3mm}\noindent {\bf Remark.} If $A\otimes M\to M$ is a representation of a graded algebra in a graded vector space $M$ then we define $M^\vee \otimes A\to M^\vee$ by
$$\langle m, f\cdot a\rangle = (-1)^{deg\, a} \langle a\cdot m, f\rangle.$$
This is a right representation of $A$: the \emph{dual right representation}. In particular the holonomy representation dualizes to a right representation of $\overline{UL_V}$ in $H(\land Z)^\vee$.

\vspace{3mm} In the special case that $\land W = \land V\otimes \land Z$ is itself a minimal Sullivan algebra,   the homotopy Lie algebras form a short exact sequence
$$0 \leftarrow L_V\leftarrow L_W\leftarrow L_Z\leftarrow 0$$
identifying $L_Z$ as an ideal in $L_W$. In particular, for $\beta\geq \alpha\in {\mathcal J}$, $L_{Z_\alpha}$ is a finite dimensional ideal in $L_{W_\beta}$, and the right adjoint representation of $L_{W_\beta}$ in $L_{Z_\alpha}$ is nilpotent. As above, these representations extend to a right representation of $\overline{UL_W}$ in $L_{Z_\alpha}$, and passing to inverse limits defines a representation of $\overline{UL_W}$ in $L_Z$, extending the right adjoint representation of $L_W$ dual to the adjoint representation.

\vspace{3mm}\noindent {\bf Definition.} This representation, denoted $ad_R$,  is the \emph{right adjoint representation} of $\overline{UL_W}$ in $L_Z$.

\vspace{3mm} Again recall from (\ref{l13}) that the Hurewicz map $\xi^\vee$ induces a linear map $L_Z/L_Z^{(2)} \to H(\land Z)^\vee$ of degree 1.

\begin{Prop}\label{prop7} Suppose $\land W= \land V\otimes \land Z$ decomposes a minimal Sullivan algebra as a $\Lambda$-extension of a minimal Sullivan algebra, $\land V$. Then \begin{enumerate}
\item[(i)] Each $L_Z^{(r)}$ is a sub $\overline{UL_W}$-module.
\item[(ii)] The surjection $L_W\to L_V$ extends to a surjection $\pi : \overline{UL_W}\to \overline{UL_V}$.
\item[(iii)] The quotient right representations of $\overline{UL_W}$ in $L_Z^{(r)}/L_Z^{(r+1)}$ factor over $\pi$ to yield right representations of $\overline{UL_V}$.
\item[(iv)] The Hurewicz map $\xi^\vee : L_Z/L_Z^{(2)} \to H^{\geq 1}(\land Z)^\vee$ is a morphism of $\overline{UL_V}$-modules.
\end{enumerate}\end{Prop}

\noindent {\sl proof.} (i) This is immediate, because if $\beta\geq \alpha $ then each $L^r_{Z_\alpha}$ is a $\widehat{UL_{W_\beta}}$-module.

(ii) This is immediate because $\pi$ is the inverse limit of the surjections, 
$$\pi_{r,\alpha} : UL_{W_\alpha}/I_\alpha^r \to UL_{V\alpha}/J_\alpha^r,$$
between finite dimensional vector spaces, where $I_\alpha$ and $J_\alpha$ respectively denote the augmentation ideals in $UL_{W_\alpha}$ and $UL_{V_\alpha}$.

(iii) This follows because in each $\land V_\alpha\otimes \land Z_\alpha$, the adjoint representation of $L_{Z_\alpha}$ is zero in $L^r_{Z_\alpha}/L^{r+1}_{Z_\alpha}$.

(iv) Recall that division by $\land^{\geq 2}Z$ is the linear map $\land^{\geq 1}Z\to Z$ which induces $\xi: H^{\geq 1}(\land Z)\to Z$.   Because $\land V\otimes \land Z$  is itself a minimal Sullivan algebra, the derivations $\theta_i$ of (\ref{l15})  preserve $\land^{\geq 1}Z$. Define $\widehat{\theta}_i : Z\to Z$ by requiring 
$$\xi\circ \theta_i = \widehat{\theta}_i\circ \xi.$$
Then for $x\in L_W$, $y\in L_Z$ and $z\in Z$ we have, where $v_i$ is a basis of $V$,
$$\langle d_1(1\otimes z), sx,sy\rangle = \sum \, \langle v_i\otimes \widehat{\theta}_iz, sx,sy\rangle = \sum (-1)^{(deg\, x+1)(deg\, y+1)}\langle v_i, s\overline{x}\rangle \langle \widehat{\theta}_iz, sy\rangle,$$
where $\overline{x}$ is the image of $x$ in $L_V$. 

Now suppose $z= \xi\Phi$ for some $\overline{d}$-cycle $\Phi\in \land Z$. Then $\widehat{\theta}_iz= \xi \theta_i \Phi$. It follows that
$$\langle d_1(1\otimes z), sx,sy\rangle =- (-1)^{(deg\, y+1)(deg\, x+1)}\langle \xi\overline{\theta}\Phi, sy\rangle.$$
On the other hand, by definition
$$\langle d_1(1\otimes z), sx,sy\rangle = (-1)^{deg\, y+1} \langle z, [x,y]\rangle = (-1)^{deg\, y+1} \langle z, ad\,x(y)\rangle.$$
A straightforward computation now gives that
\begin{eqnarray}\label{l16}
\xi^\vee \circ ad_R(\overline{x}) = \overline{\theta}^\vee(\overline{x}) \circ \xi^\vee.
\end{eqnarray}

Finally (iv) follows from (\ref{l16}) applied to the sub $\Lambda$-extension of the form $\land V_\alpha \otimes \land Z_\alpha$ in which dim$\, V_\alpha \oplus Z_\alpha <\infty$, together with a standard limit argument.
 
 \hfill$\square$

 \vspace{3mm}\noindent {\bf Example.} \emph{Acyclic Closures}
 
 Suppose $(\land V\otimes \land U,d)$ is the acyclic closure of a minimal Sullivan algebra, $\land V$. This is the $\Lambda$-extension of $\land V$ satisfying $\varepsilon : \land V\otimes \land U \stackrel{\simeq}{\to} \mathbb Q$, where $\varepsilon$ is the augmentation sending $V,U\to 0$. It has the important property that $d : U\to \land^+V\otimes \land U$, so that the quotient differential in $\land U$ is zero. Thus in this case the holonomy representation is a representation of $\overline{UL_V}$ in $\land U$. This then dualizes to a right representation of $\overline{UL_V}$ in $(\land U)^\vee$. 
 
 Now let $\varepsilon_U: \land U\to \mathbb Q$ be the augmentation vanishing on $\land^+U$. Then we have
 
\begin{Prop}  An isomorphism of right $\overline{UL_V}$-modules
 $$\overline{UL_V} \stackrel{\cong}{\longrightarrow} (\land U)^\vee$$
 is given by $a\mapsto \varepsilon_U\cdot a$.
 \end{Prop}
 
 \vspace{3mm}\noindent {\sl proof:}  Suppose $V_\alpha\subset V$ is a finite dimensional subspace for which $dV_\alpha \subset \land V_\alpha$. Then the inclusion $V_\alpha \hookrightarrow V$ extends to an inclusion of acyclic closures $\land V_\alpha \otimes \land U_\alpha \hookrightarrow \land V\otimes \land U$ which maps $U_\alpha$ to a subspace of $U$. Now because dim$\, V_\alpha<\infty$, (\cite[Theorem 6.1]{RHTII}) asserts that
 $$\overline{UL_{V_\alpha}} \stackrel{\cong}{\to} (\land U_\alpha)^\vee.$$
 Passing to inverse limits gives the isomorphism of the Proposition, because $\land U= \varinjlim_\alpha \land U_\alpha$ and so $(\land U)^\vee = \varprojlim_\alpha (\land U_\alpha)^\vee$. \hfill$\square$

 \section{The Sullivan rationalization, $\ell_{\mathbb Q}$}
 
 In this section all spaces are   connected and based CW complexes or simplicial sets. All maps and homotopies preserve base points, and "homotopy" is denoted "$\sim$". In particular the augmentation  $\land V\to \mathbb Q$  in a Sullivan algebra defines a base point in $\langle \land V\rangle$ and the maps $\widetilde{\varphi} : X\to \langle \land V\rangle$ of \S 1 are base point preserving with respect to any base point of $X$.
 
 Here we define homotopy localization,   exhibit Sullivan's construction as an example, and examine its properties from that perspective. Homotopy localizations are defined via   the category $\mathcal C$ described next, when there is no ambiguity, \emph{we may use the same symbol to denote a map and its homotopy class}.
 \begin{enumerate}
 \item[$\bullet$] The objects in ${\mathcal C}$ are the families $(X):= \{X_\alpha\}$ of connected based spaces, together with the assignment  for each pair $X_\alpha, X_\beta\in (X)$ of a single homotopy class of homotopy equivalences,
 $$\omega_{\alpha, \beta} : X_\beta \stackrel{\cong}{\longrightarrow} X_\alpha.$$
 These are required to satisfy $\omega_{\alpha,\beta}\circ \omega_{\beta, \gamma}\sim \omega_{\alpha, \gamma}$ and $\omega_{\sigma, \sigma} \sim \mbox{id}_{X_\sigma}$.
 \item[$\bullet$] A morphism $(g): (X)\to (Y)$ in ${\mathcal C}$ is the assignment for each pair $X_\alpha\in (X), Y_\beta\in (Y)$ of a homotopy class of maps,
 $$g_{\alpha,\beta} : X_\beta\to Y_\alpha,$$
 satisfying $g_{\alpha,\beta}\circ \omega_{\beta,\gamma}^X \sim \omega^Y_{\alpha,\delta}\circ g_{\delta, \gamma}$ for all $\alpha, \beta,\gamma, \delta$.
 \end{enumerate}
 
 \vspace{3mm}\noindent {\bf Remark.} 
 Given objects $(X)$ and $(Y)$ in ${\mathcal C}$, a single map $g_{\beta,\alpha}: X_\alpha \to Y_\beta$ extends uniquely to a morphism $(g)$ in ${\mathcal C}$. The map $g_{\beta,\alpha}$ is called a \emph{representative} for $(g)$. A morphism $(g)$ is a \emph{homotopy equivalence} if some $g_{\beta, \alpha}$ (equivalently all $g_{\beta, \alpha}$) is a homotopy equivalence.
 
 Moreover, the maps $\omega_{\alpha, \beta}$ provide a consistent identification of the spaces $X_\alpha\in (X)$ as a single homotopy type, and the maps $g_{\alpha, \beta}$ provide a consistent identification of $(g)$ as the homotopy class of a map $(X)\to (Y)$.

 \vspace{3mm}
 In particular, if $(X)\in {\mathcal C}$ then the homotopy equivalences $\omega_{\alpha, \beta}$ identify the groups $\pi_i(X_\alpha)$ with a single group $\pi_i(X)$ and the graded algebras $H(X_\alpha)$ with a single graded algebra $H(X)$. This then defines functors
 $$\pi_i : {\mathcal C}\leadsto {\mathcal G} \hspace{5mm}\mbox{and } H: {\mathcal C} \leadsto {\mathcal A}$$
 to the category of groups and graded algebras. In particular, the various constructions of the classifying space of a group $G$ define a natural transformation
 $B : {\mathcal G} \longrightarrow {\mathcal C}$ 
 from the category of groups, characterized by
  $\pi_1\circ B = id$,  and $ \pi_i\circ B = 0$, $i\geq 2.$ 
 
 On the other hand, the homotopy category of spaces and homotopy classes of maps is a subcategory of ${\mathcal C}$ via a canonical functor,
 $$h : {\mathcal H}\leadsto {\mathcal C},$$
 defined as follows:
  \begin{enumerate}
 \item[$\bullet$] $hX$ contains the single space $X$,
 \item[$\bullet$] ${\mathcal C}(hX,hY)$ is the set of homotopy classes of maps from $X$ to $Y$.
 \end{enumerate}
 
 \vspace{3mm}\noindent {\bf Definition.} A \emph{homotopy localization} is a functor $\ell : {\mathcal H}\leadsto {\mathcal C}$ together with a natural transformation $t : h\to \ell$.
  Thus if $\ell X= \{X_\sigma, (\omega_{\sigma, \tau})\}$ then $t(X)$ is a family of homotopy classes of maps $f_\sigma: X\to X_\sigma$ satisfying $\omega_{\sigma, \tau}\circ f_\tau \sim f_\sigma$. 
 
 \vspace{3mm} For the definition of Sullivan rationalization it is convenient to extend the definition of Sullivan models to quasi-isomorphisms $\land W\to A_{PL}(X)$ from $\Lambda$-algebras. Because our spaces are connected such a Sullivan model always decomposes as a tensor product $(\land V,d)\otimes \land (U\oplus dU)$ in which $\land V$ is a minimal Sullivan algebra and $d : U\stackrel{\cong}{\rightarrow} dU$. In particular, this implies that the inclusion $\land V\to \land W$ induces a homotopy equivalence $\langle \land V\rangle \leftarrow \langle \land W\rangle$.
 
  \vspace{3mm}\noindent {\bf Definition.} The \emph{Sullivan rationalization} is the natural transformation $t_{\mathbb Q}: h\to \ell_{\mathbb Q}$, defined next.
   
 \begin{enumerate}
 \item[$\bullet$] $\ell_{\mathbb Q}(X)= \{\langle \land V_\sigma\rangle\}$, indexed by the Sullivan models $\varphi_\sigma :  \land V_\sigma \stackrel{\simeq}{\to} A_{PL}(X)$,    together with the   homotopy equivalences $\langle \varphi_{\tau, \sigma}\rangle : \langle \land V_\sigma\rangle \to \langle \land V_\tau\rangle$ determined by the unique class of morphisms $\varphi_{\tau, \sigma}: \land V_\sigma \stackrel{\simeq}{\to} \land V_\tau$ satisfying $\varphi_\tau \circ \varphi_{\tau, \sigma} \sim \varphi_\sigma$.  
 \item[$\bullet$] If $f : X\to Y$ then $\ell_{\mathbb Q}f= \{ \langle \psi_{\tau, \sigma}\rangle : \langle \land V_\sigma\rangle \to \langle \land W_\tau\rangle\}$, where the $\psi_{\tau, \sigma}$ are the Sullivan representatives of $f$.
 \item[$\bullet$] The natural transformation $t_{\mathbb Q}: h\to \ell_{\mathbb Q}$ assigns for each $X$ the homotopy classes of the maps $\widetilde{\varphi_\sigma}: X\to \langle \land V_\sigma\rangle$ defined in \S 1.
 \end{enumerate}

 \vspace{3mm}\noindent {\bf Remark.} If $\ell_{\mathbb Q}(X)= \ell_{\mathbb Q}(Y)$ then $X$ and $Y$ have the same Sullivan models. Moreover, if $f: X\to Y$,  the single homotopy class determined by $\ell_{\mathbb Q}(f)$ is a homotopy equivalence if and only if a Sullivan representative $\varphi : \land V\to \land W$ of $f$ is a quasi-isomorphism.
   
 \vspace{3mm} Our principal objective in this section (Proposition 8) is to characterize Sullivan rationalization, without reference to Sullivan models,   in terms of
 \begin{enumerate}
 \item[$\bullet$] its behaviour on \emph{elementary spaces}, defined next, and
 \item[$\bullet$] its behaviour with respect to \emph{inverse systems}.
 \end{enumerate}

 \vspace{3mm}\noindent {\bf Definition.} An \emph{elementary space} is a nilpotent space, $X$, such that each $\pi_1(X)^r/\pi_1(X)^{r+1}$ is a finite dimensional rational vector space, and $\oplus_{i\geq 2} \pi_i(X)$ is also a finite dimensional rational vector space.
 
 \vspace{3mm} As for inverse systems, we consider only those whose index set is directed set ${\mathcal S}$ satisfying
\begin{eqnarray}
\label{r17}
\mbox{For each $\gamma\in {\mathcal S}$ there is an $n$ such that if $\gamma_1<\dots <\gamma_r= \gamma$, then $r\leq n$.}
\end{eqnarray}
 For a directed set ${\mathcal S}$ satisfying (\ref{r17}), denote by $n(\gamma)$, $\gamma \in \mathcal S$, the least $n$ for which (\ref{r17}) holds. Then ${\mathcal S} = \cup_n {\mathcal S}(n)$, where
 $${\mathcal S}(n) = \{ \gamma \in {\mathcal S}\,\vert\, n(\gamma) \leq n\}.$$
Where there is no ambiguity we denote $\displaystyle\varprojlim_{\gamma \in {\mathcal S}}$ and $\displaystyle\varinjlim_{\gamma\in {\mathcal S}}$ by $\varprojlim$ and $\varinjlim$.

\begin{Prop} \label{prop8} Sullivan rationalization satisfies
\begin{enumerate}
\item[(i)] If $X$ is an elementary space,   then $t_{\mathbb Q}(X): X\to \ell_{\mathbb Q}(X)$ is a homotopy equivalence.   
\item[(ii)] If $f : X\to Y$ represents $t_{\mathbb Q}(X): X\to \ell_{\mathbb Q}X$ then there is a map $g = \{g_\sigma\}: Y\to \varprojlim  Y_\sigma$ in which the $Y_\sigma$ are elementary spaces and the composite
$$\varinjlim  H(Y_\sigma) \to H(Y) \to H(X)$$
is an isomorphism.
 \item[(iii)] If a map $g= \{g_\sigma\} : X\to  \varprojlim  X_\sigma$ induces an isomorphism 
 $$H(X)\stackrel{\cong}{\longleftarrow} \varinjlim  H(X_\sigma)$$
 then
 $$\pi_*(\ell_{\mathbb Q}X)\stackrel{\cong}{\longrightarrow} \varprojlim  \pi_*(\ell_{\mathbb Q}X_\sigma).$$
 \item[(iv)] If $\widehat{t} : h\to \widehat{\ell}$ is any homotopy localization satisfying conditions (i)-(iii), then for any space $X$, $\ell_{\mathbb Q}(X)$ and $\widehat{\ell}(X)$ have the same homotopy type.  
  \end{enumerate}
 \end{Prop}

 \noindent {\bf Lemma 1.} The following conditions on a connected space $X$ are equivalent
\begin{enumerate}
\item[(i)] $X$ is an elementary space.
\item[(ii)] The minimal Sullivan model of $X$ has the form $\land V$ with dim$\, V^i<\infty$, and the map $X\to \langle \land V\rangle$ is a homotopy equivalence.
\end{enumerate}

\vspace{3mm}\noindent {\sl proof:} 

(i) $\Rightarrow$ (ii)  Denote $\pi_1(X)$ by $G$ and suppose $G^{n+1}= \{0\}$. If $G$ is abelian then $BG = S^1_{\mathbb Q}\times \dots \times S^1_{\mathbb Q}= \langle \land V\rangle$, where $\land V\to A_{PL}(BG)$ is the minimal Sullivan model. In general (\cite[Theorem 5.1]{RHTII}) gives a commutative diagram
$$\xymatrix{
A_{PL}(B\, G/G^2) \ar[r] & A_{PL}(BG) \ar[r] & A_{PL}(BG^2)\\
\land W\ar[u]\ar[r] & \land W\otimes \land Z \ar[u]\ar[r] & \land Z\ar[u]}$$
in which the vertical arrows are Sullivan models. Apply $\langle \,\rangle$ to get a map of fibrations
$$\xymatrix{
B\, G/G^2\ar[d] & BG\ar[l]\ar[d] & BG^2\ar[l]\ar[d]\\
\langle \land W\rangle   & \langle \land W\otimes \land Z\rangle \ar[l] & \langle \land Z\rangle.  \ar[l]}$$
By induction on $n$, and the case $n=1$, the left and right arrows are homotopy equivalences. Thus so is the central arrow.

Finally, the same argument applied to the fibration
$$BG\leftarrow \widehat{X}\leftarrow \widetilde{X}$$
where $\widehat{X}\simeq X$ and $\widetilde{X}$ is the universal cover, shows that if $\land V$ is the minimal Sullivan model of $X$ then $X\to \langle \land V\rangle$ is a homotopy equivalence. 

(ii) $\Rightarrow $ (i)  If (ii) holds then $\pi_i(X)= (V^i)^\vee$ and so (i) follows at once from \cite[Chap. 2]{RHTII}. 

\hfill$\square$

\vspace{3mm}\noindent {\bf Lemma 2.} If $g_{\sigma, \tau}: X_\sigma\leftarrow X_\tau$, $\sigma\leq \tau \in {\mathcal S}$, defines an inverse system of maps, then there is a family of commutative diagrams
$$\xymatrix{
A_{PL}(X_\sigma) \ar[rr]^{A_{PL}(g_{\sigma, \tau})} && A_{PL}(X_\tau)\\
\land V_\sigma \ar[u]_{\varphi_\sigma}^\simeq \ar[rr]_{\psi_{\sigma, \tau}} && \land V_\tau \ar[u]_\simeq ^{\varphi_\tau} , && \sigma \leq \tau \in {\mathcal S}}$$
in which
\begin{enumerate}
\item[(i)] the $\psi_{\sigma, \tau}$ form an inductive system of morphisms of $\Lambda$-algebras;
\item[(ii)] each $\psi_{\sigma, \tau}$ restricts to an inclusion $V_\sigma \to V_\tau$;
\item[(iii)] each $\varphi_\sigma$ is a quasi-isomorphism.
\end{enumerate}

\vspace{3mm}\noindent {\sl proof:} The construction is by induction. If $\tau \in {\mathcal S}(1)$ we let $\varphi_\tau : \land V_\tau \to A_{PL}(Y_\tau)$ be any Sullivan model. Suppose the construction is accomplished for $\tau \in {\mathcal S}(n)$ and let $\tau \in {\mathcal S}(n+1)/{\mathcal S}(n)$. Then set $\overline{V}_\tau = \varinjlim_{\sigma <\tau} \land V_\sigma$, so that 
$$\varinjlim_{\sigma <\tau}\varphi_\sigma : \land \overline{V}_\tau \to \varinjlim_{\sigma<\tau} A_{PL}(V_\sigma)$$
is a quasi-isomorphism from a $\Lambda$-algebra satisfying $H^0(\land \overline{V}_\tau)= \mathbb Q$.

Because $Y_\tau$ is connected it follows from (\cite[Theorem 3.1]{RHTII}) that the composite
$$\land \overline{V}_\tau \to \varinjlim_{\sigma<\tau} A_{PL}(Y_\sigma) \to A_{PL}(Y_\tau)$$
extends to a quasi-isomorphism $\land \overline{V}_\tau\otimes \land Z \stackrel{\simeq}{\to} A_{PL}(X_\tau)$ from a $\Lambda$-algebra. This is the desired quasi-isomorphism $\varphi_\tau : \land V_\tau \to A_{PL}(Y_\tau)$, and   this closes the induction.
\hfill$\square$

\vspace{3mm}\noindent {\bf Lemma 3.} Suppose $\land V = \varinjlim  \land V_\sigma$ is the direct limit of an inductive system of $\Lambda$-algebras in which the morphisms map $V_\sigma \to V_\tau$. If each $H^0(\land V_\sigma)= \mathbb Q$ then $\pi_*\langle \land V\rangle \stackrel{\cong}{\longrightarrow} \varprojlim \pi_*\langle \land V_\sigma\rangle$. 

\vspace{3mm}\noindent {\sl proof:} Define differentials $d_0$ in $V$ and $V_\sigma$ by setting $d_0v$ to be the component of $dv$ in $V$( or in $V_\sigma$). Then in the minimal models $\land Z$ and $\land Z_\sigma$ for $\land V$ and $\land V_\sigma$ we have isomorphisms which identify
$$Z_\sigma \cong H(V_\sigma, d_0) \hspace{5mm}\mbox{and } Z\cong H(V, d_0) = \varinjlim  H(V_\sigma, d_0).$$
These isomorphisms identify $\pi_*\langle \land V_\sigma\rangle= H(V_\sigma, d_0)^\vee$ and $\pi_*\langle \land V\rangle = H(V, d_0)^\vee$. Since dualizing converts direct limits to inverse limits, it follows that
$$\pi_*\langle \land V\rangle = \left[ \varinjlim  H(V_\sigma, d_0)\right]^\vee = \varprojlim  \pi_*\langle \land V_\sigma\rangle.$$
\hfill$\square$

\vspace{3mm}\noindent {\sl proof of Proposition \ref{prop8}:} 
(i) This is immediate from Lemma 1.

\vspace{2mm} (ii) 
We may assume $f= \widetilde{\varphi}: X\to \langle \land V\rangle$, where $\varphi: \land V \to A_{PL}(X)$ is a Sullivan model. This identifies $H(f): H(Y)\to H(X)$ with $H(\widetilde{\varphi})$. Now as described 
 in \S 3, $\land V = \varinjlim_\sigma \land V_\sigma$, where the $\land V_\sigma$ are sub Sullivan algebras and dim$\, V_\sigma<\infty$. Thus by Lemma 1, each $Y_\sigma = \langle \land V_\sigma\rangle$ is a Sullivan space. Moreover, since 
  $\langle \, \rangle$ converts direct limits to inverse limits, $Y= \varprojlim Y_\sigma$.   Finally, by   \cite[Theorem 5.4]{RHTII}   adjoint to $\mbox{id}_{\land V_\sigma}$ is the quasi-isomorphism $\varphi_\sigma : \land V_\sigma \stackrel{\simeq}{\rightarrow} A_{PL}\langle \land V_\sigma\rangle$. Since by the definition $\varphi = \varinjlim_\sigma \varphi_\sigma$ it follows from (3) that the composite
  $$\varinjlim H(Y_\sigma) \to H(Y)\to H(X)$$
  is the isomorphism $H(\varphi)$.

\vspace{2mm} (iii) Lemma 2 provides an inductive system of quasi-isomorphisms $\varphi_\sigma : \land V_\sigma \to A_{PL}(X_\sigma)$ in which
$$\varphi : \land V = \varinjlim \land V_\sigma \stackrel{\simeq}{\to} \varinjlim_\sigma A_{PL}(X_\sigma) \stackrel{\simeq}{\to} A_{PL}(X)$$
  exhibits $\widetilde{\varphi} : X\to \langle \land V\rangle$ as a representative of $t_{\mathbb Q}(X).$ Moreover, the quasi-isomorphisms $\varphi_\sigma$ in Lemma 2 identify $\pi_*(\ell_{\mathbb Q}X_\sigma)$ with $\pi_*\langle \land V_\sigma\rangle$ and $\pi_*(\ell_{\mathbb Q}X)$ with $\pi_*\langle \land V\rangle$. This identifies $\pi_*(\ell_{\mathbb Q}X)\to \varprojlim \pi_*(\ell_{\mathbb Q}X_\sigma)$ with $\pi_*\langle \land V\rangle\to \varprojlim_\gamma \pi_*\langle \land V_\sigma\rangle$. This is an isomorphism by Lemma 3.

\vspace{2mm} (iv) Suppose $\widehat{t}: h\to \widehat{l}$ is a second functor and natural transformation satisfying (i)-(iii). Then, for a space $X$, let
$$\xymatrix{X\ar[rr]^{f=\widehat{t}(X)} && Y=\widehat{X} \ar[rr]^{\{g_\sigma\}} && \varprojlim  Y_\sigma}$$
be the maps provided by (ii). Then Lemma 2 provides a  morphism
$$\xymatrix{\psi: \land V = \varinjlim_\sigma  \land V_\sigma \ar[rr]_\simeq^{\varinjlim_\sigma \varphi_\sigma} && \varinjlim_\sigma A_{PL}(Y_\sigma)\ar[r] & A_{PL}(Y)}$$
from a $\Lambda$-algebra $\land V$. Moreover, by hypothesis,  
$$  A_{PL}(f)\circ \psi : \land V\stackrel{\simeq}{\to} A_{PL}(X) $$ is a Sullivan model.

On the other hand,   we have the commutative diagram,
$$
\xymatrix{ \pi_*(Y) \ar[d]^{\pi_*\widetilde{\psi}} \ar[rr]^{\{\pi_*g_\sigma\}} && \varprojlim \pi_*(Y_\sigma)\ar[d]^{\{\pi_*\widetilde{\varphi}_\sigma\}}\\
\pi_*\langle \land V\rangle \ar[rr] && \varprojlim \pi_*\langle \land V_\sigma\rangle.}
$$
Now $Y = \widehat{X}$ and, since each $Y_\sigma$ is an elementary space, by (i) there is a natural identification $\pi_*(Y_\sigma) = \pi_*(\widehat{Y_\sigma})$. Thus by (iii), the upper horizontal arrow is an isomorphism.

On the other hand, since each $Y_\sigma$ is an elementary space, Lemma 1 asserts that each $\widetilde{\varphi}_\sigma$ is a homotopy equivalence. Thus the right hand vertical arrow is an isomorphism. Finally, Lemma 3 asserts that the lower horizontal arrow is an isomorphism. Therefore, $\pi_*(\widetilde{\psi})$ is an isomorphism and $\widetilde{\psi}$ is a homotopy equivalence. Since $A_{PL}(f)\circ\psi : \land V\stackrel{\simeq}{\to} A_{PL}(X)$ it follows that $\langle \land V\rangle \in \ell_{\mathbb Q}(X)$ and 
$$\widetilde{\psi} : \ell_{\mathbb Q}X \stackrel{\simeq}{\to} \widehat{X}.$$

\hfill$\square$

\vspace{3mm} Recall now from \cite{BK} that the Bousfield-Kan completion $\mathbb Q_{\infty}(X)$ can be expressed as the inverse limit
$$X\to \varprojlim_k\mathbb Q_k(X) = \mathbb Q_\infty(X)$$
of a tower of fibrations in which $\varinjlim_k H(\mathbb Q_k(X)) \stackrel{\cong}{\longrightarrow} H(X)=2^k $. Moreover, Proposition \ref{prop8} is analogous to \cite[Tower Lemma 6.2]{BK}. Here we provide a proof of the following result of Bousfield:

\begin{Prop}\label{prop9} For any connected space, the map $X\to \ell_{\mathbb Q}(X)$ factors up to homotopy as
$$X\to \mathbb Q_{\infty}(X) \to \ell_{\mathbb Q}(X).$$\end{Prop}

\vspace{3mm}\noindent{\sl Proof:}   Apply Lemma 2   to $X\to\mathbb Q_\infty(X)= \varprojlim_k \mathbb Q_k(X)$ to obtain a quasi-isomorphism,
$$\varphi: \land V= \varinjlim \land V_k \to A_{PL}(\mathbb Q_\infty(X)) \to A_{PL}(X),$$
which then gives the map $\mathbb Q_\infty (X)\to \langle \land V\rangle = \ell_{\mathbb Q}(X)$.

 \hfill$\square$

\vspace{3mm}\noindent {\bf Corollary.} If $H(\mathbb Q_\infty (X)) \stackrel{\cong}{\to} H(X)$ - ($X$ is $\mathbb Q$-good in the terminology of \cite{BK}) - then $\mathbb Q_\infty (X)_{\mathbb Q} \stackrel{\cong}{\to}   X_{\mathbb Q}$.

\vspace{3mm} Finally, Proposition \ref{prop8} provides a proof of a major theorem of Bousfield and Gugenheim:

\begin{Prop} \label{prop10}(\cite[Theorem 12.2]{BG})   If $X$ is a connected space and $H(X)$ is a graded vector space of finite type, then $\mathbb Q_\infty(X)\to X_{\mathbb Q}$ is a homotopy equivalence.\end{Prop}

 \noindent {\sl proof:} As observed in the third example in \S 3, in a minimal Sullivan model $\land V$ for $X$, $V$ admits an increasing filtration $V[0]\subset \dots \subset V[n]\subset \dots$ in which $d: V[n+1]\to \land V[n]$ and each $V[n]$ is finite dimensional. 
Set $X[n]= \langle \land V[n]\rangle$. Then $X_{\mathbb Q}$ is the inverse limit of the tower of fibrations
$$X[0]\leftarrow \dots \leftarrow X[n]\leftarrow$$
But this tower satisfies   conditions of (\cite[Chapter 5]{BK}) which identify the inverse limit as $\mathbb Q_\infty(X)$. \hfill$\square$

 \vspace{3mm}\noindent {\bf Remark.}  Suppose that for some $n$,   card $H_n(X;\mathbb Q)$ is an infinite cardinal $k= $ card$\, \oplus_i H_i(X_{\mathbb Q})$. Then Bousfield has shown \cite{2a} that card $\mathbb Q_\infty (X) \leq 2^k$. On the other hand, if $\land V$ is a minimal Sullivan model of $X$, then card$\,V=$ card $\land V= $ card$\, H(X)= 2^k$. Then
 $$\mbox{card}\, \pi_*(X_{\mathbb Q})= \mbox{card}\, V^\vee = 2^{2^k} > \mbox{card}\, \mathbb Q_\infty (X).$$
 In particular, the hypothesis of finite type in the result of Bousfield and Gugenheim is necessary.

 \vspace{5mm}\noindent Institut de Math\'ematique et de Physique, Universit\'e Catholique de Louvain, 2, Chemin du cyclotron, 1348 Louvain-La-Neuve, Belgium, yves.felix@uclouvain.be

 \vspace{1mm}\noindent Department of Mathematics, Mathematics  Building, University of Maryland, College Park, MD 20742, United States, shalper@umd.edu
 
 \end{document}